\documentclass[final,onefignum,onetabnum]{siamart190516}

\newsiamremark{remark}{Remark}
\newsiamremark{hypothesis}{Hypothesis}
\crefname{hypothesis}{Hypothesis}{Hypotheses}
\newsiamthm{claim}{Claim}

\headers{The AAA\MakeLowercase{trig} algorithm}{P. J. Baddoo}

\title{The AAA\MakeLowercase{trig} algorithm for Rational Approximation of Periodic Functions
}

\author{Peter J. Baddoo\thanks{Department of Mathematics, Imperial College London, South Kensington Campus 
  (\email{p.baddoo@imperial.ac.uk}, \url{https://www.baddoo.co.uk}).}}

\usepackage{amsopn}
\DeclareMathOperator{\diag}{diag}

\usepackage{lipsum}
\usepackage{amsfonts}
\usepackage{graphicx}
\usepackage{epstopdf}
\usepackage{algorithmicx,algpseudocode}
\ifpdf
  \DeclareGraphicsExtensions{.eps,.pdf,.png,.jpg}
\else
  \DeclareGraphicsExtensions{.eps}
\fi
\usepackage{todo}

\usepackage{mathtools}
\usepackage{amssymb}
\usepackage{graphicx}
\usepackage{mwe}
\usepackage{subcaption}
\renewcommand{\i}{\textnormal{i}}
\renewcommand{\d}{\textnormal{d}}
\newcommand{\e}{\textnormal{e}}

\newcommand{\bC}{\boldsymbol{C}}
\newcommand{\bA}{\boldsymbol{A}}
\newcommand{\bw}{\boldsymbol{w}}

\newcommand{\bS}{\boldsymbol{S}}
\renewcommand{\Re}{\mathbb{R}\textnormal{e}}
\renewcommand{\Im}{\mathbb{I}\textnormal{m}}

\DeclarePairedDelimiter\floor{\lfloor}{\rfloor}
\newcommand{\cst}{\textnormal{cst}}
\usepackage{amsmath}
\DeclareMathOperator*{\argmax}{arg\,max}

\newcommand{\norm}[1]{\left\lVert#1\right\rVert}

\usepackage{tikz,pgfplots,pstool}

\usetikzlibrary{calc, 
				decorations.markings, 
				decorations.text,
				decorations.pathreplacing,
				shapes.multipart,
				positioning,
				arrows.meta,
				arrows,
				matrix,
				patterns,
				shadows,
				hobby}

\tikzset{
    side by side/.style 2 args={
        line width=2pt,
        #1,
        postaction={
            clip,postaction={draw,#2}
        }
    },
    circle node/.style={
        circle,
        draw,
        fill=white,
        minimum size=1.3cm
    }
}

\tikzset{test/.style={
    postaction={
        decorate,
        decoration={
           markings,
            mark=at position \pgfdecoratedpathlength-0.5pt with {\arrow[blue,line width=#1] {}; },
            mark=between positions 0 and \pgfdecoratedpathlength-0pt step 0.5pt with {
                \pgfmathsetmacro\myval{multiply(divide(
                    \pgfkeysvalueof{/pgf/decoration/mark info/distance from start}, \pgfdecoratedpathlength),100)};
                \pgfsetfillcolor{white!\myval!black};
                \pgfpathcircle{\pgfpointorigin}{#1};
                \pgfusepath{fill};}
}}}}

\tikzset{
    mark position/.style args={#1(#2)}{
        postaction={
            decorate,
            decoration={
                markings,
                mark=at position #1 with \coordinate (#2);
            }
        }
    }
}

\def\myShiftUp#1{\raisebox{1ex}}
\def\myShiftDown#1{\raisebox{-2.5ex}}

\def\myPathTextLeftAbove#1#2#3#4{
	\draw [LaTeX-,thick,postaction={decorate,
		decoration={
			raise=1ex,
			text along path,
			text align={center},
			text={%
				|\color{black}| {#1} }
		}
	}
	] #2 to [bend left=#4] #3;
}

\def\myPathTextLeftBelow#1#2#3#4{
	\draw [LaTeX-,thick,postaction={decorate,
		decoration={
			raise=-2ex,
			text along path,
			text align={center},
			text={%
				|\color{black}| {#1} }
		}
	}
	] #2 to [bend left=#4] #3;
}

\usetikzlibrary{decorations.text,fit,chains,calc,shapes.geometric,intersections}

\newlength{\myww}
\newlength{\myhh}

\usepackage{amsfonts}

\tikzset{->-/.style={decoration={
  markings,
  mark=at position #1 with {\arrow{latex}}},postaction={decorate}}}

\newlength{\fheight}
\newlength{\fwidth}

\ifpdf
\hypersetup{
  pdftitle={The AAAtrig algorithm for rational approximation of periodic functions},
  pdfauthor={P. J. Baddoo}
}
\fi

\begin{document}
\maketitle

\begin{abstract}
We present an extension of the AAA (adaptive Antoulas--Anderson) algorithm for periodic functions, called `AAAtrig'. 
The algorithm uses the key steps of AAA approximation by
(i) representing the  approximant in (trigonometric) barycentric form and
(ii) selecting the support points greedily.
Accordingly, AAAtrig inherits all the favourable characteristics of AAA
and is thus extremely flexible and robust, being able to consider quite general sets of sample points in the complex plane.
We consider a range of applications with particular emphasis on solving Laplace's equation in periodic domains
and compressing periodic conformal maps.
These results reproduce the tapered exponential clustering effect observed in other recent studies.
The algorithm is implemented in Chebfun.
\end{abstract}

\begin{keywords}
	rational functions, trigonometric rational functions, AAA algorithm, harmonic functions, conformal mapping
\end{keywords}

\begin{AMS}
	41A20, 65D15
\end{AMS}
 \section{Introduction} \label{Sec:introduction}
%

Given a collection of sample points and data values from a function $f$, 
can we stably construct a rational function that approximates $f$?
Can the approximation be improved by exploiting knowledge of the structure of $f$?
The former question has recently been addressed with the AAA (adaptive Antoulas--Anderson) algorithm \cite{Nakatsukasa2018};
the present work addresses the latter question when $f$ is known to be \textit{periodic}.
Periodic functions are ubiquitous in every area of science and there is significant scope for application of rational approximants.
For example, PDEs are commonly formulated in periodic domains, which reflects the fact that many 
relevant physical domains are, or can be modelled as, spatially periodic.
Periodic rational approximants can thus be used to design effective spectral methods \cite{Tee2006}
or tackle elliptic PDEs directly using the recently proposed `lightning solvers' \cite{Gopal2019a}.
In other areas, discovering periodicities in the time series of a signal 
is important for analysis and diagnostics in engineering applications \cite{Bendat2012}.
Periodic rational function approximations can compress such signals and also uncover underlying structure \cite{Wilber2020}.


The AAA algorithm was introduced as a robust, fast and flexible method for generating rational approximants from data.
Given its effectiveness, it is surprisingly simple: 
the original article presented a \textsc{Matlab} implementation in 40 lines of code.
Various extensions of AAA have been developed for minimax approximation \cite{Nakatsukasa2019},
matrix-valued functions \cite{Gosea2020}, surrogate functions \cite{Elsworth2019} 
and parametric dynamical systems \cite{Rodriguez2020}.
In the present work we present an enhancement of AAA for approximating periodic (trigonometric) functions, called `AAAtrig'.
A central idea of AAA algorithms is to use the barycentric representation of 
rational functions \cite{Berrut2004,Berrut2014}.
As opposed to representing the rational approximant as a quotient of polynomials $p(z)/q(z)$, which can be ill-conditioned,
AAA represents the approximant as a ratio of partial fractions:
\begin{align}
	r(z) = 
	\left. \sum_{j = 1}^m \frac{\alpha_j}{z - z_j} 
	\middle/
	\sum_{j = 1}^m \frac{\beta_j}{z - z_j} ,
	\right.
	\label{Eq:baryInt}
\end{align}
where $\alpha_j$, $\beta_j$ and $z_j$ are constants chosen to optimize the approximation.
This barycentric representation is better conditioned than both the aforementioned polynomial quotient 
and the partial fractions formula used by vector fitting \cite{Gustavsen1999}.
The trigonometric analogue of the barycentric formula follows quite naturally as
\begin{align}
	r(z) = 
	\left. \sum_{j = 1}^m \alpha_j \cst\left(\frac{z - z_j}{2}\right)
	\middle/
	\sum_{j = 1}^m \beta_j \cst\left(\frac{z - z_j}{2}\right),
	\right.
	\label{Eq:baryInt2}
\end{align}
where $\cst$ is a trigonometric function defined in section \ref{Sec:algorithm} \cite{Henrici1979}.
Rational approximation can seem an intimidating nonlinear problem due to the unknown support points $\{z_j\}$.
The other main step in the AAA algorithm addresses this challenge by selecting the support points $\{z_j\}$ 
from the sample points greedily, which has the additional effect of mitigating exponential instabilities.
Again, this step generalises naturally to the trigonometric case.
Thus, AAAtrig is identical to AAA except the polynomial basis functions in \eqref{Eq:baryInt}
have been replaced with trigonometric basis functions \eqref{Eq:baryInt2}.
It follows that AAAtrig inherits all the favourable properties of AAA: it is fast,
it can be applied to a range of domains, and it appears to give an optimal distribution of poles and zeros.
It is unusual for numerical Froissart doublets (very near poles and zeros) to appear, and when they do,
they can be removed using a similar `cleanup' procedure to that of AAA.

Rational functions approximation is particularly useful when approximating functions with singularities.
This is perhaps illustrated most clearly by the canonical result of Newman \cite{Newman1964}
who showed that $|x|$ can be approximated with root-exponential accuracy by rationals, 
whereas polynomials achieve merely algebraic accuracy.
Another appealing feature of rational functions is that they perform well on unbounded domains;
the other canonical problem of approximating $\e^{-x}$ on $[0,\infty)$ serves at a good example \cite{Cody1969}.
These observations have recently motivated the 
application of rational functions to solve PDEs with singularities \cite{Gopal2019a}.
These lightning solvers, which belong to the class of Methods of Fundamental Solutions (MFS) \cite{Fairweather1998},
have been successfully applied to solve the Laplace, Helmholtz \cite{Gopal2019c} and biharmonic equations.
Once the solution is computed, it can be compressed with the AAA algorithm to generate a compact representation
of the solution that is extremely rapid to evaluate.
Typically these compressed solutions can be evaluated hundreds of thousands of points in a second.
In this paper we generalise these ideas to solve PDEs in periodic domains and compress the solutions using AAAtrig.
We apply this theory to calculate the potential flow past a periodic array of obstacles.

Since these lightning solvers can be used to solve Laplace's equation,
they can also be used to calculate conformal maps \cite{Gopal2019b,Trefethen2020Conf}.
Indeed, the motivation of the present author to develop AAAtrig was to numerically compress conformal maps found using the recent
periodic Schwarz--Christoffel formula \cite{Baddoo2019c}.
This compression means that the forwards and backwards maps can be computed with typical speed-ups of between
10--1,000 times compared to evaluating the actual Schwarz--Christoffel formula, which requires evaluating numerous integrals.
We also show how the AAA framework can be used to compress periodic and multi-valued conformal maps.

The power of rational approximation for these tasks stems from the exponential clustering of poles and zeros near singularities,
as is apparent from Newman's explicit construction \cite{Newman1964}.
It has recently been proposed that the optimal distribution of these poles and zeros follows a tapered distribution
where the logarithm of their density linearly decays to zero near the singularity \cite{Trefethen2020Singularity}.
Our findings support this hypothesis 
and we show that AAAtrig selects tapered poles and zeros in our conformal mapping application.




The remainder of the paper is arranged as follows.
In section \ref{Sec:algorithm} we introduce the basic ideas of trigonometric interpolation and then
explicate the AAAtrig algorithm.
In section \ref{Sec:examples} we consider a number of applications of AAAtrig including comparisons with
AAA and FFT-based interpolation, 
the removal of Froissart doublets,
solutions of Laplace's equation in periodic domains, 
and compression of conformal maps.
We conclude in section \ref{Sec:conclusions} with a summary and discussion of future work.
The algorithm is available in Chebfun \cite{Chebfun} as \texttt{aaatrig}.

 \section{The AAAtrig algorithm} \label{Sec:algorithm}
In this section we will present the AAAtrig algorithm.
We will first review trigonometric interpolation and justify our form of rational approximant.
We then outline the key steps of AAAtrig that determine the support points and weights.
Following this, we explain how to find the poles and zeros once the support points
and weights have been determined.

\subsection{The form of the approximant}
We begin by outlining the general theory of trigonometric polynomial and rational approximation.
A thorough summary of the underlying mathematics of computations with trigonometric functions is available in \cite{Wright2015}.
A trigonometric polynomial of degree $m$ with period $2 \pi$ is a function of the \mbox{form \cite{Henrici3}}
\begin{align}
	r(z) &= \sum_{k = -m}^m c_k \e^{\i k z}.
	\label{Eq:trigPoly}
\end{align}
Alternatively, $r$ may be expressed in terms of sines and cosines as
\begin{align}
	r(z) &= \frac{1}{2} a_0 + \sum_{k=1}^m a_k \cos(k z) + b_k \sin(k z).
\end{align}
Our dependent variable is represented here as $z$, as opposed to the more traditional $\theta$.
This choice of notation is used to emphasize the flexibility of AAAtrig:
the sample points are not tied to a particular domain such as the real interval $[0,2 \pi)$.

For $m$ function values $f_j$ at sample points $z_j$, the Lagrangian representation of trigonometric polynomial 
that interpolates $f$ is
\begin{align}
	r(z) = \sum_{j=1}^m f_j a_j l(z) \cst\left( \frac{z - z_j}{2}  \right),
	\label{Eq:trigInterp}
\end{align}
where $l$ is the trigonometric node polynomial
\begin{align}
	l(z) = \prod_{j=1}^m \sin \left( \frac{z - z_j}{2} \right) ,
	\label{Eq:trigNodePoly}
\end{align}
and the coefficients are 
\begin{align}
	a_j = \prod_{\substack{k \neq j\\k = 0} }^m \csc \left( \frac{z_k - z_j}{2} \right),
\end{align}
where
\begin{align}
	\cst(z) &=
	\begin{cases}
		\csc(z) & \textnormal{if $m$ is odd,}\\
		\cot(z) + c & \textnormal{if $m$ is even.}
	\end{cases}
	\label{Eq:cstDef}
\end{align}
The distinction between odd and even numbers of sample points stems from the 
odd \mbox{($2 m +1$)} number of unknowns in \eqref{Eq:trigPoly}.
As such, the interpolation problem for an even number of interpolation points is under-determined.
The problem is usually closed by specifying the form of the highest frequency of oscillation
by requiring that the coefficient of $\sin(m z)$ vanishes \cite{Henrici1979},
and the trigonometric polynomial is `balanced' \cite{Berrut1984}.
In this case, that would mean $c = \cot( \sum {z_j} /2)$, but we take $c=0$ in our calculations.
The reason for this choice is that we wish to consider arbitrary sample data;
obviously there will be a problem if, for example, $\sum z_j = 2 \pi k$ for integer $k$.

Dividing \eqref{Eq:trigInterp} by the corresponding expression for $f \equiv 1$ yields an alternative expression for $r$:
\begin{align}
	r(z) = \frac{\sum_{j=1}^m f_j a_j l(z) \cst\left( \frac{z - z_j}{2}  \right)}
{\sum_{j=1}^m a_j l(z) \cst\left( \frac{z - z_j}{2}  \right)}
= 
 \frac{\sum_{j=1}^m f_j a_j \cst\left( \frac{z - z_j}{2}  \right)}
 {\sum_{j=1}^m a_j \cst\left( \frac{z - z_j}{2}  \right)}. \label{Eq:bary1}
\end{align}
The above expression is  the \emph{trigonometric barycentric formula} for $r(z)$ \cite{Berrut1984}.
The barycentric form has several advantages over the Lagrange formula \eqref{Eq:trigInterp},
as outlined in \cite{Berrut2004}.
First, evaluating the Lagrangian form requires for each $z$ requires $\mathcal{O}(m^2)$ floating point operations (flops), 
whereas evaluating the  barycentric form \eqref{Eq:bary1} requires $\mathcal{O}(m)$ flops 
once the coefficients $w_j$ have been computed.
Second, including a new data pair ($f_{m+1}, z_{m+1}$) requires every computation to be performed again
in the Lagrangian formulation; 
the barycentric form can be updated in $\mathcal{O}(m)$ flops.
 Higham \cite{Higham2004} has proved that the (polynomial) barycentric formula is forwards stable, 
 as long as the interpolating points have a small Lebesgue constant.
Third, Berrut proved that under certain conditions, the trigonometric barycentric formula \eqref{Eq:bary1}
is guaranteed to be well-condition for interpolation on the unit circle \cite{Berrut1988}.
One drawback of the barycentric form is that locations of the zeros of $r$ are now obscured.

The barycentric form \eqref{Eq:bary1} interpolates $f$ regardless of the choice of coefficients $a_j$.
Accordingly, we may replace the $a_j$ with arbitrary weights $w_j$ to obtain
\begin{align}
	r(z) = 
\left.
	\sum_{j=1}^m f_j w_j \cst\left( \frac{z - z_j}{2}  \right) \middle/
	 \sum_{j=1}^m w_j \cst\left( \frac{z - z_j}{2}  \right). 
 \right. \label{Eq:approximant}
\end{align}
The function $r$ is a trigonometric polynomial in the special case $w_j = a_j$.
For any other choices of $w_j$, the function $r$ represents a trigonometric rational function, i.e.,
the ratio of two trigonometric polynomials.
We write the numerator and denominator of \eqref{Eq:approximant} as $n$ and $d$ respectively,
so that $r(z) = n(z)/d(z)$.
The barycentric form is related to the traditional form of a rational function (i.e. $r(z) = p(z)/q(z)$)
by multiplication of the trigonometric node polynomial $l(z)$ \eqref{Eq:trigNodePoly}.


%
%
Before introducing the AAAtrig algorithm, we educe the correspondence between the trigonometric and polynomial
barycentric forms.
The barycentric form for a (non-trigonometric) rational interpolating $f_j$ at $z_j$ is \cite{Berrut2004}
\begin{align}
	r(z) &=	
	\left.
\sum_{j=1}^m \frac{f_j w_j}{z - z_j} 
		\middle/
		\sum_{j=1}^m \frac{w_j}{z - z_j}. 
		\right.
		\label{Eq:polyBary}
\end{align}
This is the form of ansatz used in the original AAA algorithm.
The connection between the barycentric forms for trigonometric and polynomial interpolants 
can be made more transparent by expanding $\cst$ from \eqref{Eq:cstDef} in its Mittag--Lefler representation.
In particular, substituting
\begin{align}
	\csc(z) = \sum_{k = -\infty}^\infty \dfrac{(-1)^k}{z - k \pi}, \qquad \qquad
	\cot(z) = \sum_{k = -\infty}^\infty \dfrac{1}{z - k \pi}. 
\end{align}
into the rational approximant \eqref{Eq:approximant} yields
\begin{align}
	r(z) &= \left.
		\sum_{k = -\infty}^\infty \chi_{m}^k \left( \sum_{j = 1}^m \frac{f_j w_j}{z - (z_j - 2 \pi k)} \right)
		\middle/
		\sum_{k = -\infty}^\infty \chi_{m}^k \left( \sum_{j = 1}^m \frac{f_j w_j}{z - (z_j - 2 \pi k)} \right)
		\right.
		,
		\label{Eq:expansion}
\end{align}
where $\chi_{m} = (-1)^{(m+1)}$.
As such, it can be seen that a AAAtrig approximant \eqref{Eq:polyBary} 
can be viewed as a AAA approximant with periodic support points.

\subsection{The algorithm}
In this section we present the core AAAtrig algorithm.
The guiding principles follow directly from the original AAA formulation and only minor changes are required.
It's worth discussing the choice of basis function $\cst$.
Our problem is distinct from that considered by Henrici \cite{Henrici1979} and Berrut \cite{Berrut1984}
since our order of approximation varies with each iteration, whereas they consider interpolants of fixed order.
Accordingly, the definition of $\cst$ in \eqref{Eq:cstDef} is not so relevant, as the definition
changes with each iteration.
In the present algorithm, the form of the approximant is specified to
be `odd' or `even' by the user.
Odd approximants correspond to $\cst \equiv \csc$ whereas even approximants correspond to $\cst \equiv \cot$.
Numerical experiments indicate that $\cot$ or $\csc$ are equally good choices
as far as the speed of convergence is concerned%
\footnote{Alternating expressions were also considered but there was no discernible
	improvement in performance and the implementation was much more complicated.}%
	.

We consider a sample set $Z \in \mathbb{C}$ of finite size $M$.
We additionally consider a function $f(z)$ that is defined for all $z \in Z$. 
Since the underlying function $f$ is periodic, we project the sample points onto the strip
$0\leq \Re[z] <2 \pi$ by
\begin{align}
	z_j \mapsto z_j - 2 \pi \Bigl\lfloor\Re\left[\frac{z_j}{2\pi}\right]\Bigr\rfloor
\end{align}
for $j = 1, \dots, M$.

We now come to the first essential step of the AAA framework. 
At step $m$, the next support point is selected `greedily':
the support point $z_m$ is chosen to minimise the error between the function values and 
the approximant evaluated at the remaining support points.
In other words,
\begin{align}
	z_j = \argmax_{z \in Z^{(m-1)}} \norm{f(z) - \frac{n(z)}{d(z)}}
\end{align}
where $Z^{(m-1)}$ is the set of sample points with all the previous support points removed, i.e.
\begin{align}
	Z^{(m-1)} = Z\setminus \{z_1, \cdots, z_{m-1} \}.
\end{align}
Having found the sample point with the maximum residual error and adding it to the set of support points,
we now seek to calculate the weights $w_1, \cdots, w_m$ that solve the least-squares problem
\begin{align}
	\min_{z \in Z^{(m)}} \norm{d f - n}_2, \qquad \qquad \norm{\bw}_2 = 1.
	\label{Eq:minW}
\end{align}
As such, we require $m\leq M/2$ so that $Z^{(m)}$ contains at least $m$ points.
We now detail the linear algebra involved in solving this least squares problem.

We define the column vectors of the weights
\begin{align*}
	\bw = \left[ w_1, \cdots, w_m \right]^{T}.
\end{align*}
The minimisation problem \eqref{Eq:minW} can equivalently be expressed as finding $\bw$ such that 
\begin{align}
	\min \norm{\bA^{(m)} \bw}_2, \qquad \qquad \norm{\bw}_2 = 1,
	\label{Eq:leastSquares}
\end{align}
where
\begin{align*}
	\bA^{(m)} = 
	\begin{bmatrix}
	\left(F_1^{(m)} - f_1\right) \cst\left( \dfrac{Z_1^{(m)} - z_1}{2} \right) & 
	\cdots & 
	\left(F_1^{(m)} - f_m\right) \cst\left( \dfrac{Z_1^{(m)} - z_m}{2} \right)  \\[5ex]
	\vdots & \ddots & \vdots \\[5ex]
	\left(F_{M-m}^{(m)} - f_1\right) \cst\left( \dfrac{Z_{M-m}^{(m)} - z_1}{2} \right) &
	\cdots & 
	\left(F_{M-m}^{(m)} - f_m\right) \cst\left( \dfrac{Z_{M-m}^{(m)} - z_m}{2} \right) 
	\end{bmatrix}
\end{align*}
is the $(M-m) \times m$ trigonometric analogue of the Loewner matrix.
Since $M$ can be very large, it is expedient to exploit sparsity and write
$\bA^{(m)} = \bS_F \bC - \bC \bS_f$
where 
\begin{align}
	\bS_F = \diag \left(F_1^{(m)}, \cdots, F_{M-m}^{(m)}\right), \qquad 
	\bS_f = \diag \left(f_1, \cdots, f_{M-m}\right),
\end{align}
and
$\bC$ is the trigonometric $(M-m)\times m$ Cauchy matrix
\begin{align*}
\bC = 
\begin{bmatrix}
\cst\left( \dfrac{Z_1^{(m)} - z_1}{2} \right) & 
\cdots & 
 \cst\left( \dfrac{Z_1^{(m)} - z_m}{2} \right)  \\[5ex]
\vdots & \ddots & \vdots \\[5ex]
\cst\left( \dfrac{Z_{M-m}^{(m)} - z_1}{2} \right) &
\cdots & 
\cst\left( \dfrac{Z_{M-m}^{(m)} - z_m}{2} \right) 
\end{bmatrix}.
\end{align*}
%

This concludes the description of the core AAAtrig algorithm.
Readers familiar with the original AAA formulation \cite{Nakatsukasa2018} will note that 
AAAtrig is almost identical to AAA.
The only substantial difference arises in the definition of $\bC$, where we have replaced the polynomial
basis with a trigonometric basis. 
\subsection{Far-field behaviour}
\label{Sec:farField}
In some applications it is useful to specify the behaviour of the rational approximant at infinity.
In other words, we wish to specify (finite) $f_\infty^\pm$ where
\begin{align}
	f_\infty^\pm  = \lim_{z \rightarrow \pm \i \infty} r(z).
\end{align}
Odd approximants are of the form \eqref{Eq:trigApproxOdd} and $r$ takes different values at
$\pm \i \infty$:
\begin{align}
	f_\infty^\pm &= \frac{\sum_{j=1}^m f_j w_j \e^{\mp z_j/2}}{\sum_{j=1}^m w_j \e^{\mp z_j/2}}.
	\label{Eq:ffOdd}
\end{align}
This constraint can be included in the least squares problem \eqref{Eq:leastSquares} by
simply appending the matrix $\bA^{(m)}$ with the $2 \times m$ matrix
\begin{align*}
	\begin{bmatrix}
		\left(f_\infty^+ - f_1\right) \e^{-\i z_1/2} & \cdots & \left(f_\infty^+ - f_m\right) \e^{-\i z_m/2}\\
		\left(f_\infty^- - f_1\right) \e^{+\i z_1/2} & \cdots & \left(f_\infty^- - f_m\right) \e^{+\i z_m/2}
	\end{bmatrix}.
\end{align*}
Even approximants are of the form \eqref{Eq:trigApproxEven} and
$r$ takes the same value as $z \rightarrow \pm \i \infty$ so that
\begin{align}
	f_\infty^\pm = f_\infty &= \frac{\sum_{j=1}^m f_j w_j}{\sum_{j=1}^m w_j}.
	\label{Eq:ffEven}
\end{align}
Again, including this condition requires only a minor modification to the least-squares problem \eqref{Eq:leastSquares};
simply appending $\bA^{(m)}$ with the row vector
\begin{align*}
	\begin{bmatrix}
		f_\infty - f_1 & \cdots & f_\infty - f_m
	\end{bmatrix}
\end{align*}
results in the desired system.

\subsection{Finding the poles and zeros} \label{Sec:polesZeros}
The core AAAtrig algorithm does not provide the poles and zeros of the rational function
--- it only gives the support points and weights. 
The poles and zeros and zeros are found in the following way.
In what follows, we will make use of the exponential representations of $\csc$ and $\cot$:
\begin{align}
	\csc(z) = \frac{2 \i \e^{\i z}}{\e^{2 \i z} - 1},
	\qquad
	\cot(z)	= \i \left( 1 + \frac{2}{\e^{2\i z} - 1} \right).
	\label{Eq:cstExp}
\end{align}
\subsubsection{Odd approximants}
In this case, the trigonometric rational approximant is
\begin{align}
	r(z) = \left. \sum_{j=1}^m f_j w_j \csc \left( \frac{z - z_j}{2}  \right) \middle/
		\sum_{j=1}^m  w_j \csc \left( \frac{z - z_j}{2}  \right). \right.
		\label{Eq:trigApproxOdd}
\end{align}
We can manipulate the trigonometric approximant \eqref{Eq:trigApproxOdd} into the polynomial
approximant \eqref{Eq:polyBary} by a suitable transformation.
In particular, inserting the substitutions 
\begin{align}
\hat{z} = \e^{\i z}, \qquad \hat{z}_j = \e^{\i z_j}, \qquad \hat{w}_j = w_j \e^{\i z_j/2},
\label{Eq:oddSub}
\end{align}
into \eqref{Eq:trigApproxOdd} yields
\begin{align}
	r(z) \triangleq R(\hat{z}) = \left. \sum_{j=1}^m \frac{f_j \hat{w}_j}{ \hat{z} - \hat{z}_j}  \middle/
	\sum_{j=1}^m  \frac{\hat{w}_j}{\hat{z} - \hat{z}_j} \right. ,
\end{align}
where we have used the identity (\ref{Eq:cstExp}a).
%
%
The poles and zeros may now be found in the same way as the original AAA algorithm \cite{Nakatsukasa2018}
using a method originally proposed by Klein \cite{Klein2012}.
For example, the zeros of $R$ are the eigenvalues of the $(m+1)\times(m+1)$ generalised eigenvalue problem
\begin{align}
	\begin{pmatrix}
	0 & f_1 \hat{w}_1 & \cdots & f_m \hat{w}_m\\
	1 &\hat{z}_1&&\\
	\vdots &&\ddots&\\
	1 &&&\hat{z}_m
	\end{pmatrix} \boldsymbol{v} = \lambda
	\begin{pmatrix}
	0 & & &\\
	& 1 & &\\
	& & \ddots &\\
	& & & 1
	\end{pmatrix}
	\boldsymbol{v} 
	\label{Eq:oddEig}
\end{align}
with eigenvectors
\begin{align*}
	\boldsymbol{v} = \left[1 , 1/(\lambda - \hat{z}_1),\cdots, 1/(\lambda - \hat{z}_m) \right]^{T}.
\end{align*}
At least two of the $m+1$ eigenvalues of the system \eqref{Eq:oddEig} are infinite and these are discarded
since a type $(m,m)$ approximant can have, at most, $m-1$ zeros.
The poles are the eigenvalue of an identical problem with $f_j \hat{w}_j$ replaced with $\hat{w}_j$.
Once the zeros and poles of $R$ have been found, 
the zeros and poles of $r$ can be found by inverting the substitution \eqref{Eq:oddSub}.
\subsubsection{Even approximants}
Finding the poles and zeros for barycentric trigonometric rational functions follows in a similar manner
upon use of a suitable transformation.
In this case, the rational barycentric form is
\begin{align}
	r(z) = \left. \sum_{j=1}^m f_j w_j \cot \left( \frac{z - z_j}{2}  \right) \middle/
		\sum_{j=1}^m  w_j \cot \left( \frac{z - z_j}{2}  \right). \right.
		\label{Eq:trigApproxEven}
\end{align}
Applying
\begin{align}
\tilde{z} = \tan\left(\frac{z}{2}\right) \label{Eq:trans}
\end{align}
converts the trigonometric barycentric form \eqref{Eq:trigApproxEven} to
\begin{align}
	r(z) \triangleq R(\tilde{z})= \left. \left({\sum_{j=1}^m \frac{f_j \tilde{w}_j}{\tilde{z} - \tilde{z}_j}} + c_n \right)\middle/\left({\sum_{j=1}^m \frac{\tilde{w}_j}{\tilde{z} - \tilde{z}_j}} + c_d \right)\right., \label{Eq:baryTrans}
\end{align}
where the support points and weights are now
\begin{align}
 \tilde{z}_j = \tan\left(\frac{z_j}{2}\right), \qquad \qquad \tilde{w}_j = w_j(1+\tilde{z}_j^2),  \label{Eq:transSupportWeights}
\end{align}
and the constants in the numerator and denominator of $R$ are respectively given by
\begin{align}
c_n = \sum_{j=1}^m  \frac{f_j \tilde{w}_j \tilde{z}_j}{1+\tilde{z}_j^2}, \qquad c_d = \sum_{j=1}^m \frac{\tilde{w}_j\tilde{z}_j}{1+\tilde{z}_j^2}. \label{Eq:constDef}
\end{align}
The special case where $z_j=\pi$ is a support point is dealt with in section \ref{Sec:specialSupport}.

Note that \eqref{Eq:baryTrans} is in the original barycentric form \eqref{Eq:polyBary} 
with additional constants in the numerator and denominator. 
Accordingly, the poles and zeros of $R$ can be found by solving a generalised eigenvalue problem \cite{Klein2012}. 
For example, it is straightforward to verify that the eigenvalues of the generalised eigenvalue problem
\begin{align}
	\begin{pmatrix}
	c_d & f_1 \tilde{w}_1 & \cdots & f_m \tilde{w}_m\\
	1 &\tilde{z}_1&&\\
	\vdots &&\ddots&\\
	1 &&&\tilde{z}_m
	\end{pmatrix} \boldsymbol{v} = \lambda
	\begin{pmatrix}
	0 & & &\\
	& 1 & &\\
	& & \ddots &\\
	& & & 1
	\end{pmatrix}
	\boldsymbol{v} \label{Eq:genEig}
\end{align}
are the poles of $R$, and the eigenvectors are given by
\begin{align*}
\boldsymbol{v} = \left[1 , 1/(\lambda - \tilde{z}_1),\cdots, 1/(\lambda - \tilde{z}_m) \right]^{T}.
\end{align*}
The zeros of $r$ are then given by inverting the transformation \eqref{Eq:trans}.
The poles of $r$ can be found by solving a similar generalised eigenvalue problem where 
$c_d$ is replaced by $c_n$ and the $\tilde{w}_j$ are replaced by $f_j \tilde{w}_j$ in \eqref{Eq:genEig}.

\subsubsection*{Special case: $\boldsymbol{\pi}$ is a support point} \label{Sec:specialSupport}
\newcommand{\pSum}{\sum_{j=2}^m}
In the above, the symbolic manipulations break down when the approximant is even and
$z_j = \pi$ is a support point for some $j$.
Assuming this to be the case, we take $j=1$ so that $z_1 = \pi$.
The approximant \eqref{Eq:trigApproxEven} can therefore be rewritten as 
\begin{align*}
r(z) =  \Bigg( -f_1 w_1 \tan\left(\frac{z}{2} \right) +
& \pSum f_j w_j \cot \left(\frac{z - z_j}{2}\right) \Bigg) \Bigg/\\
 \Bigg(  - w_1 \tan\left(\frac{z}{2} \right) + 
&\pSum w_j \cot \left(\frac{z - z_j}{2}\right) \Bigg) .
\end{align*}
Now applying the transformation \eqref{Eq:trans} yields
\begin{align}
	r(z) \triangleq R(\tilde{z})= \left. \left( \frac{c_n}{\tilde{z}} - f_1 w_1 +
	\pSum \frac{f_j \tilde{w}_j}{\tilde{z}(\tilde{z} - \tilde{z}_j)} \right)\middle/
	\left(\frac{c_d}{\tilde{z}} - w_1 + 
	\pSum \frac{\tilde{w}_j}{\tilde{z}(\tilde{z} - \tilde{z}_j)} \right)\right. \label{Eq:specialR}
\end{align}
where $c_n$ and $c_d$ are defined in \eqref{Eq:constDef} with the exception that the
first term is now excluded in the summation. 
Finding the poles and zeros of \eqref{Eq:specialR} can now be formulated as a generalized eigenvalue problem.
For example, the zeros of \eqref{Eq:specialR} are the eigenvalues of the problem
\begin{align}
	\begin{pmatrix}
	-f_1 w_1 & c_d & f_2 \tilde{w}_2, &\cdots & f_m \tilde{w}_m\\
		1 &0 &&&\\
	0 & 1& \tilde{z}_2 &\\
	\vdots & \vdots&& \ddots \\
	0 & 1 &&&\tilde{z}_m
	\end{pmatrix} \boldsymbol{v} = \lambda
	\begin{pmatrix}
		0 & & & &\\
		  & 1 & & &\\
		  & & \ddots & &\\
	& & & \ddots &\\
	& & & & 1
	\end{pmatrix}
	\boldsymbol{v}
\end{align}
which has eigenvectors
\begin{align*}
	\boldsymbol{v} = \left[1 , \frac{1}{\lambda}, 
	\frac{1}{\lambda(\lambda - \tilde{z}_2)},\cdots, \frac{1}{\lambda(\lambda - \tilde{z}_m)} \right]^{T}.
\end{align*}
The poles can be found in a similar manner.
%

Once the poles and zeros have been found, the residues at each pole can be computed with the standard quotient formula.
Then the approximant can alternatively be expressed in its partial fraction representation as
\begin{align}
	r(z) = \sum_{k=1}^{m-1} q_k \cot\left(\frac{z - p_k}{2}\right)  + c
	\label{Eq:partial}
\end{align}
where $\{p_j\}$ are the poles, $q_j$ is the residue at $p_j$, and $c$ is a constant.
Comparison with the far-field values in section \ref{Sec:farField} shows that,
for odd approximants, the residues and constant must satisfy \eqref{Eq:ffOdd}
\begin{align}
	c \mp \i \sum_{k=1}^{m-1} q_k = 
	 \frac{\sum_{j=1}^m f_j w_j \e^{\mp z_j/2}}{\sum_{j=1}^m w_j \e^{\mp z_j/2}}.
\end{align}
whereas for even approximants we have \eqref{Eq:ffEven}
\begin{align}
	\sum_{k=1}^{m-1} q_k = 0, \qquad \qquad c = \frac{\sum_{j=1}^m f_j w_j}{\sum_{j=1}^m w_j}.
\end{align}
From the above we see that the even approximant is constrained such that its residues sum to zero.
In contrast, the odd $\csc$ approximant is slightly more general and should therefore usually be favoured in calculations.
Whilst the partial fraction representation \eqref{Eq:partial} is useful for clarifying the structure of the odd and even
approximants the conversion between barycentric form and partial fractions is usually unstable.
For example, to differentiate the approximant one should not differentiate \eqref{Eq:partial} symbolically
but instead use the differentiation matrices of \cite{Baltensperger2002}.


 \section{Examples} \label{Sec:examples}
We now present a range of examples and applications of AAAtrig.
%
\subsection{Comparison with AAA}
%
%
\begin{figure}
	\centering
	\begin{subfigure}{.45\textwidth}
	\setlength{\fwidth}{\linewidth}
	\setlength{\fheight}{4cm}
	\centering
%
%
\begin{tikzpicture}[%
trim axis left, trim axis right
]

\begin{axis}[%
width=0.951\fwidth,
height=\fheight,
at={(0\fwidth,0\fheight)},
scale only axis,
xmin=1,
xlabel style={font=\color{white!15!black}},
xlabel={degree, $m$},
ymin=-14,
ymax=2,
ytick={-16,-12,-8,-4,0},
yticklabels={{$10^{-16}$},{$10^{-12}$},{$10^{-8}$},{$10^{-4}$},{$10^0$}},
ylabel style={font=\color{white!15!black}},
ylabel={maximum error},
axis background/.style={fill=white},
xmajorgrids,
ymajorgrids
]
\addplot [color=red, line width=1.0pt, forget plot]
  table[row sep=crcr]{%
1	0.438564071294124\\
2	0.250526012534774\\
3	-0.0637051529725241\\
4	-1.01778609876093\\
5	-0.992788082270547\\
6	-0.949547438636753\\
7	-3.12920191528314\\
8	-3.23126817404243\\
9	-3.53057205576441\\
10	-4.29415009717034\\
11	-4.8783638687235\\
12	-5.41015384868555\\
13	-6.35716066364234\\
14	-7.21331945567916\\
15	-8.97936794134837\\
16	-8.71151412123617\\
17	-8.95220387476141\\
18	-10.8369943413363\\
19	-12.3514693982107\\
20	-12.2411739630479\\
21	-12.6052668094523\\
};
\addplot [color=blue, line width=1.0pt, forget plot]
  table[row sep=crcr]{%
1	0.438564071294126\\
2	0.219249454816715\\
3	-0.986632249746366\\
4	-1.06509556159602\\
5	-2.92542681808408\\
6	-2.82799440210079\\
7	-5.1362309012741\\
8	-5.07551050532992\\
9	-7.62208364887027\\
10	-7.6155503014008\\
11	-10.445532088392\\
12	-10.4423093789596\\
13	-13.4001229981337\\
};
\node[right, align=left]
at (axis cs:2,-8) {AAAtrig};
\node[right, align=left]
at (axis cs:12,-4) {AAA};
\end{axis}
\end{tikzpicture}%
	\caption{$f(x) = \exp(\sin(x))$}
	\label{Fig:compa}
	\end{subfigure}
	\hfill
	\begin{subfigure}{.45\textwidth}
	\setlength{\fwidth}{\linewidth}
	\setlength{\fheight}{4cm}
	\centering
%
%
\begin{tikzpicture}[%
trim axis left, trim axis right
]

\begin{axis}[%
width=0.951\fwidth,
height=\fheight,
at={(0\fwidth,0\fheight)},
scale only axis,
xmin=1,
xlabel style={font=\color{white!15!black}},
xlabel={degree, $m$},
ymin=-12,
ymax=4,
ytick={-16,-12,-8,-4,0},
yticklabels={{$10^{-16}$},{$10^{-12}$},{$10^{-8}$},{$10^{-4}$},{$10^0$}},
ylabel style={font=\color{white!15!black}},
ylabel={maximum error},
axis background/.style={fill=white},
xmajorgrids,
ymajorgrids
]
\addplot [color=red, line width=1.0pt, forget plot]
  table[row sep=crcr]{%
1	2.72504108415583\\
2	2.22920605757118\\
3	0.46983304301297\\
4	-0.865289211363525\\
5	-2.80931772850048\\
6	-4.57800085701656\\
7	-6.9313350122915\\
8	-9.19045877936042\\
9	-11.5726871424301\\
};
\addplot [color=blue, line width=1.0pt, forget plot]
  table[row sep=crcr]{%
1	2.72504108415583\\
2	3.76499620497596\\
3	2.80268301679244\\
4	3.07462113352544\\
5	3.62429703267669\\
6	2.76964901581033\\
7	2.80570570866078\\
8	3.06326922239882\\
9	2.8010508579527\\
10	2.51074177945048\\
11	2.57811741654007\\
12	2.9280901604962\\
13	3.90452238531348\\
14	2.37129799826177\\
15	2.43008989776474\\
16	2.10540112420968\\
17	1.89300178488155\\
18	1.73617864091251\\
19	1.31270484941757\\
20	1.33035876468807\\
21	1.4514298236342\\
22	1.06223330974375\\
23	1.40192577072843\\
24	0.668220433407853\\
25	0.482542918562316\\
26	0.252531552710636\\
27	-0.095536765099439\\
28	0.114207469212801\\
29	-0.51798309332041\\
30	-0.376339611738501\\
31	-1.01936522046586\\
32	-1.33761234235662\\
33	-1.15549060072844\\
34	-1.62158358080009\\
35	-1.42506437443321\\
36	-2.20929390963572\\
37	-2.25510597294719\\
38	-2.45362485946356\\
39	-2.79525277464788\\
40	-2.90613517873017\\
41	-3.24791056411964\\
42	-3.65400758064318\\
43	-3.94858796277524\\
44	-4.10769396386991\\
45	-4.36223787125644\\
46	-4.81313788852549\\
47	-4.69631702342544\\
48	-5.50418254142879\\
49	-5.44820682531328\\
50	-6.11901206895968\\
51	-6.13614817048739\\
52	-6.73375121583375\\
53	-6.90906452682576\\
54	-7.06573262283628\\
55	-7.51396976384052\\
56	-6.50295156304487\\
57	-7.67384214411273\\
58	-6.95243394738489\\
59	-7.38882864835375\\
60	-8.94670343626795\\
61	-9.43377535727018\\
62	-9.58949192021083\\
63	-10.0133832955164\\
64	-9.77727860577785\\
65	-10.6975901990315\\
};
\node[right, align=left]
at (axis cs:40,-1) {AAAtrig};
\node[right, align=left]
at (axis cs:10,-4) {AAA};
\end{axis}
\end{tikzpicture}%
	\caption{$f(x) = \exp(x)$}
	\label{Fig:compb}
	\end{subfigure}

	\caption{Comparison of maximum errors for AAAtrig (blue) and AAA (red) for periodic and non-periodic functions.
	The functions are sampled at 1,000 randomly distributed points in a rectangle in the complex plane.
Note that the maximum error here refers to the maximum error over the set of discrete sample points $Z$.
}
	\label{Fig:AAAcomp}
\end{figure}
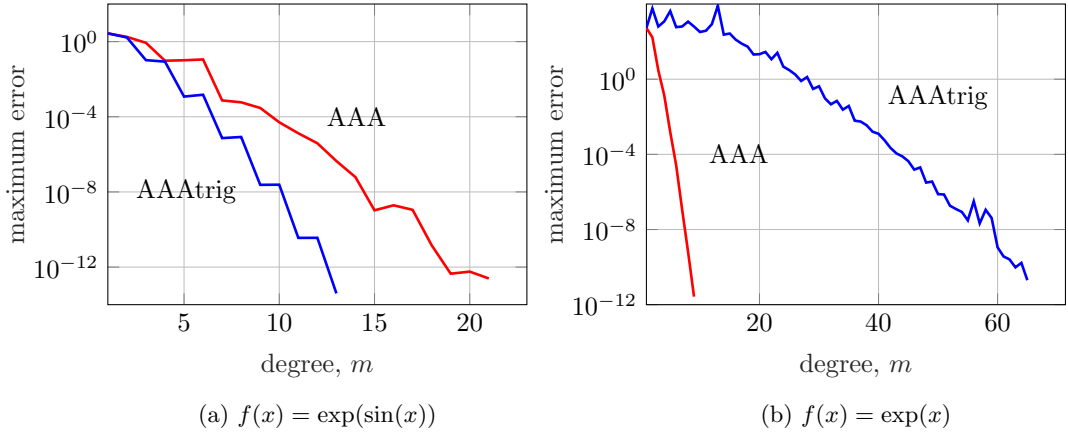
In figure \ref{Fig:AAAcomp} we compare AAAtrig to the original AAA algorithm for periodic and non-periodic functions.
The functions to be approximated, $\exp(\sin(x))$ in \ref{Fig:compa} and $\exp(x)$ in \ref{Fig:compb},
are sampled at 1,000 randomly distributed points in the rectangle $[0,2\pi]\times[-\i/2,\i/2]$.
We plot the maximum error, which is defined over the (discrete) sample set $Z$ and not a continuous domain.
As expected, AAAtrig outperforms AAA when approximating a periodic function in \ref{Fig:compa}.
Since AAAtrig uses a trigonometric barycentric approximant, it can simultaneously reduce the error at each end of the period
window.
For example, if AAAtrig chooses a support point near $0$ then the error near $2 \pi$ is automatically reduced.
Non-trigonometric AAA does not have this built-in structure and thus needs more iterations to converge.
In contrast AAA substantially outperforms AAAtrig when approximating a non-periodic function in figure \ref{Fig:compb}.
Since the approximant used by AAAtrig is constrained to be periodic, it requires many support points to resolve the 
non-periodic behaviour at the endpoints of the period window.
Note that we are using a discrete error norm so AAAtrig terminates with accuracy $10^{-11}$;
with a continuous norm we would see that AAAtrig cannot accurately represent a non-periodic function.


In these examples the improvement of AAAtrig over AAA is quite modest and the advantage of AAAtrig may not be transparent.
The differences become clearer when we wish to approximate a periodic function throughout the entire complex plane;
we encounter such examples in section \ref{Sec:harmonic}.

\subsection{Removing Froissart doublets}
Occasionally, an application of AAAtrig will return pairs of very close poles and zeros (e.g. separated by $10^{-14}$), 
known as Froissart doublets.
These anomalies can be removed using a similar `cleanup' procedure to that proposed in \cite{Nakatsukasa2018}: 
after AAAtrig has terminated, poles with residues below a certain tolerance ($10^{-13}$) 
are identified and the nearest support points are removed, followed by a final least squares problem.
This approach typically succeeds in drastically reducing the number of poles in the approximant.
\begin{figure}[tpb]
	\begin{subfigure}{\linewidth}
	\setlength{\fwidth}{\linewidth}
	\centering
	\input{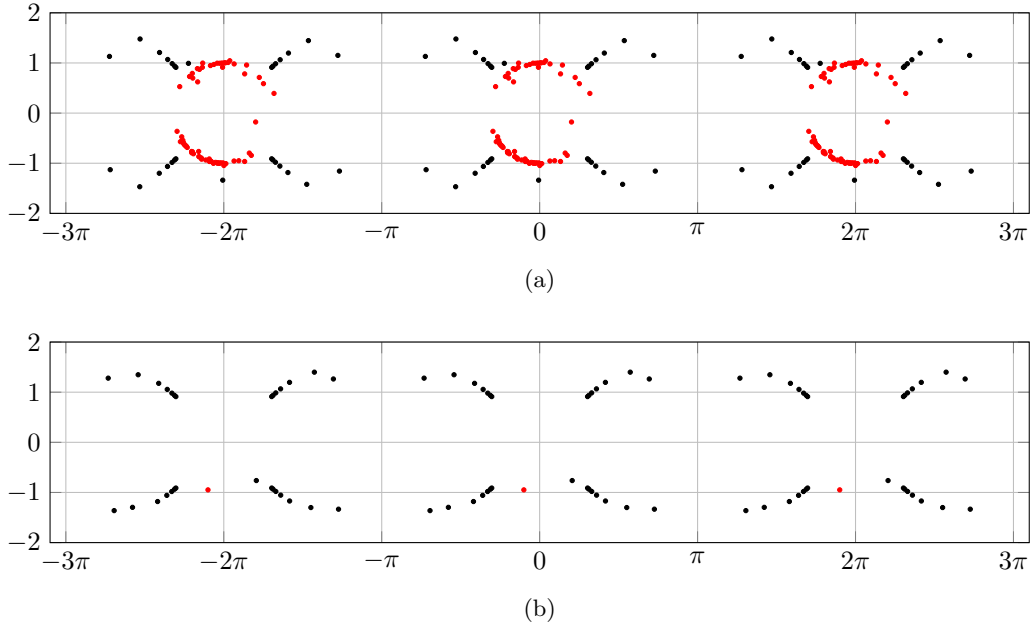}
	\caption{}
	\label{Fig:froissart1}
\end{subfigure}

	\begin{subfigure}{\linewidth}
	\setlength{\fwidth}{\linewidth}
	\centering
%
%
\begin{tikzpicture}[%
trim axis left, trim axis right
]

\begin{axis}[%
width=\fwidth,
height=0.205\fwidth,
at={(0\fwidth,0\fwidth)},
scale only axis,
xmin=-9.73893722612836,
xmax=9.73893722612836,
xtick={-9.42477796076938,-6.28318530717959,-3.14159265358979,0,3.14159265358979,6.28318530717959,9.42477796076938},
xticklabels={{$-3\pi$},{$-2 \pi$},{$-\pi$},{$0$},{$\pi$},{$2\pi$},{$3\pi$}},
ymin=-2,
ymax=2,
axis background/.style={fill=white},
xmajorgrids,
ymajorgrids
]
\addplot [color=black, draw=none, mark size=0.8pt, mark=*, mark options={solid, black}, forget plot]
  table[row sep=crcr]{%
-11.9193146527897	-0.762442511011447\\
-9.37863411292989	-3.35762815714549\\
-8.4642309501499	-1.36184349990019\\
-10.2841558741147	-1.33342821770257\\
-8.09874448307538	-1.29675583684991\\
-10.8338457478791	-1.30088751086998\\
-7.59965937273867	-1.18089870899489\\
-11.2569916062353	-1.1702502263537\\
-7.41632283816713	-1.05780381571385\\
-7.3183097106906	-0.981126279493322\\
-7.26338292441678	-0.935188857928875\\
-7.23633026538611	-0.911816452428138\\
-11.4342823370874	-1.05364642762823\\
-11.5310730805633	-0.979688590313277\\
-11.585947553353	-0.934817052430065\\
-11.6131687299628	-0.911768408188133\\
-8.58223613241985	1.27971571948419\\
-10.3858107487655	1.26308922168811\\
-7.98646171464509	1.34748890901673\\
-7.57790757983813	1.17625291666348\\
-7.4078247849691	1.05555306181727\\
-7.23613907925096	0.911791379340573\\
-7.2622491528199	0.934981309649073\\
-7.31489885973307	0.980308838476339\\
-9.47948850195751	3.2965909302135\\
-10.7644857829342	1.39901357493011\\
-11.2584062779484	1.19576257006325\\
-11.4365506704716	1.06570555254348\\
-11.5324171619942	0.98542703900096\\
-11.6132668415499	0.912175785971934\\
-11.5864839962079	0.937018234585473\\
};
\addplot [color=black, draw=none, mark size=0.8pt, mark=*, mark options={solid, black}, forget plot]
  table[row sep=crcr]{%
-5.63612934561014	-0.762442511011447\\
-3.0954488057503	-3.35762815714549\\
-2.18104564297031	-1.36184349990019\\
-4.00097056693511	-1.33342821770257\\
-1.81555917589579	-1.29675583684991\\
-4.55066044069954	-1.30088751086998\\
-1.31647406555909	-1.18089870899489\\
-4.9738062990557	-1.1702502263537\\
-1.13313753098754	-1.05780381571385\\
-1.03512440351101	-0.981126279493322\\
-0.980197617237189	-0.935188857928875\\
-0.953144958206519	-0.911816452428138\\
-5.15109702990778	-1.05364642762823\\
-5.24788777338372	-0.979688590313277\\
-5.30276224617341	-0.934817052430065\\
-5.32998342278321	-0.911768408188133\\
-2.29905082524027	1.27971571948419\\
-4.1026254415859	1.26308922168811\\
-1.70327640746551	1.34748890901673\\
-1.29472227265854	1.17625291666348\\
-1.12463947778951	1.05555306181727\\
-0.952953772071378	0.911791379340573\\
-0.979063845640309	0.934981309649073\\
-1.03171355255348	0.980308838476339\\
-3.19630319477792	3.2965909302135\\
-4.48130047575458	1.39901357493011\\
-4.97522097076877	1.19576257006325\\
-5.15336536329204	1.06570555254348\\
-5.24923185481461	0.98542703900096\\
-5.33008153437028	0.912175785971934\\
-5.30329868902829	0.937018234585473\\
};
\addplot [color=black, draw=none, mark size=0.8pt, mark=*, mark options={solid, black}, forget plot]
  table[row sep=crcr]{%
0.647055961569445	-0.762442511011447\\
3.18773650142929	-3.35762815714549\\
4.10213966420927	-1.36184349990019\\
2.28221474024447	-1.33342821770257\\
4.46762613128379	-1.29675583684991\\
1.73252486648005	-1.30088751086998\\
4.9667112416205	-1.18089870899489\\
1.30937900812388	-1.1702502263537\\
5.15004777619205	-1.05780381571385\\
5.24806090366857	-0.981126279493322\\
5.3029876899424	-0.935188857928875\\
5.33004034897307	-0.911816452428138\\
1.1320882772718	-1.05364642762823\\
1.03529753379586	-0.979688590313277\\
0.980423061006181	-0.934817052430065\\
0.953201884396377	-0.911768408188133\\
3.98413448193932	1.27971571948419\\
2.18055986559369	1.26308922168811\\
4.57990889971408	1.34748890901673\\
4.98846303452105	1.17625291666348\\
5.15854582939007	1.05555306181727\\
5.33023153510821	0.911791379340573\\
5.30412146153928	0.934981309649073\\
5.2514717546261	0.980308838476339\\
3.08688211240166	3.2965909302135\\
1.80188483142501	1.39901357493011\\
1.30796433641082	1.19576257006325\\
1.12981994388754	1.06570555254348\\
1.03395345236498	0.98542703900096\\
0.953103772809308	0.912175785971934\\
0.979886618151298	0.937018234585473\\
};
\addplot [color=black, draw=none, mark size=0.8pt, mark=*, mark options={solid, black}, forget plot]
  table[row sep=crcr]{%
6.93024126874903	-0.762442511011447\\
9.47092180860887	-3.35762815714549\\
10.3853249713889	-1.36184349990019\\
8.56540004742406	-1.33342821770257\\
10.7508114384634	-1.29675583684991\\
8.01571017365964	-1.30088751086998\\
11.2498965488001	-1.18089870899489\\
7.59256431530347	-1.1702502263537\\
11.4332330833716	-1.05780381571385\\
11.5312462108482	-0.981126279493322\\
11.586172997122	-0.935188857928875\\
11.6132256561527	-0.911816452428138\\
7.41527358445139	-1.05364642762823\\
7.31848284097545	-0.979688590313277\\
7.26360836818577	-0.934817052430065\\
7.23638719157596	-0.911768408188133\\
10.2673197891189	1.27971571948419\\
8.46374517277327	1.26308922168811\\
10.8630942068937	1.34748890901673\\
11.2716483417006	1.17625291666348\\
11.4417311365697	1.05555306181727\\
11.6134168422878	0.911791379340573\\
11.5873067687189	0.934981309649073\\
11.5346570618057	0.980308838476339\\
9.37006741958125	3.2965909302135\\
8.08507013860459	1.39901357493011\\
7.5911496435904	1.19576257006325\\
7.41300525106713	1.06570555254348\\
7.31713875954456	0.98542703900096\\
7.23628907998889	0.912175785971934\\
7.26307192533088	0.937018234585473\\
};
\addplot [color=black, draw=none, mark size=0.8pt, mark=*, mark options={solid, black}, forget plot]
  table[row sep=crcr]{%
13.2134265759286	-0.762442511011447\\
15.7541071157885	-3.35762815714549\\
16.6685102785684	-1.36184349990019\\
14.8485853546036	-1.33342821770257\\
17.033996745643	-1.29675583684991\\
14.2988954808392	-1.30088751086998\\
17.5330818559797	-1.18089870899489\\
13.8757496224831	-1.1702502263537\\
17.7164183905512	-1.05780381571385\\
17.8144315180277	-0.981126279493322\\
17.8693583043016	-0.935188857928875\\
17.8964109633322	-0.911816452428138\\
13.698458891631	-1.05364642762823\\
13.601668148155	-0.979688590313277\\
13.5467936753654	-0.934817052430065\\
13.5195724987555	-0.911768408188133\\
16.5505050962985	1.27971571948419\\
14.7469304799529	1.26308922168811\\
17.1462795140732	1.34748890901673\\
17.5548336488802	1.17625291666348\\
17.7249164437492	1.05555306181727\\
17.8966021494674	0.911791379340573\\
17.8704920758985	0.934981309649073\\
17.8178423689853	0.980308838476339\\
15.6532527267608	3.2965909302135\\
14.3682554457842	1.39901357493011\\
13.87433495077	1.19576257006325\\
13.6961905582467	1.06570555254348\\
13.6003240667241	0.98542703900096\\
13.5194743871685	0.912175785971934\\
13.5462572325105	0.937018234585473\\
};
\addplot [color=red, draw=none, mark size=0.8pt, mark=*, mark options={solid, red}, forget plot]
  table[row sep=crcr]{%
-6.59741985778967	-0.948350647294367\\
-0.31423455061008	-0.948350647294367\\
5.96895075656951	-0.948350647294367\\
12.2521360637491	-0.948350647294367\\
18.5353213709287	-0.948350647294367\\
};
\end{axis}
\end{tikzpicture}%
	\caption{}
	\label{Fig:froissart2}
\end{subfigure}
\caption{The locations of the poles in the periodic rational approximant for $f(x) = \log(2 + \cos(x)^4)$ 
for AAAtrig run with a tolerance of 0 without (a)
	and with (b) the cleanup procedure.
	The function is sampled at 1,000 roots of unity.
	Poles with residue of magnitude less than $10^{-13}$ are plotted in red;
	all other poles are plotted in black.
	The cleanup procedure removes all but one Froissart doublet.}
	\label{Fig:froissart}
\end{figure}

The effect of the cleanup procedure is illustrated in figure \ref{Fig:froissart}, 
which is comparable to figure 5.2 of \cite{Nakatsukasa2018}.
In this example, the function $f(x) = \log(2+ \cos(x)^4)$ is sampled at 1000 roots of unity and approximated with AAAtrig
with the tolerance set to zero.
(The default relative tolerance is $10^{-13}$ so this example is contrived but illustrative.)
AAAtrig identifies the branch points of $f$ and clusters poles along its branch cuts,
as illustrated in figure \ref{Fig:froissart}.
When the cleanup procedure is disabled (figure \ref{Fig:froissart1}), 
AAAtrig also identifies 66 poles which have very small numerical residue and thus represent Froissart doublets.
These poles are mainly clustered around the unit circle.
Conversely, figure \ref{Fig:froissart2} shows the effect of applying AAAtrig with the cleanup procedure enabled. 
The algorithm now terminates with 32 poles of which only one is a Froissart doublet.
The maximum error is $\mathcal{O}(10^{-13})$ in both cases.

\subsection{Comparison to FFT-based interpolation}
\label{Sec:fft}
%

The fast Fourier transform (FFT) can be used to compute a trigonometric interpolant
in the special case where the sample points are equally spaced.
For example, if we have $f_n = f( 2 \pi n/M)$ for $n = 1, \dots, M$ then we have the 
discrete Fourier transform pair
\begin{align}
	F_j = \frac{1}{M} \sum_{k=0}^{M-1} f_k \e^{-2 \pi \i k j/M},
	\qquad
	f_j = \sum_{k=0}^{M-1} F_k \e^{2 \pi \i k j/M},
\end{align}
which can be rapidly computed using a FFT algorithm.
As such, there is a family of trigonometric polynomials that interpolates the data
\begin{align}
	f(z) &= \sum_{k=0}^{M-1} F_k \e^{\i k z} \e^{\i m_k M z} \label{Eq:fftInt}
\end{align}
for any integers $m_k$.
This interpolant varies with $m_k$ and most choices of $m_k$ will result in an approximant that
oscillates wildly between the sample points.
This is the aliasing phenomenon.
How to choose $m_k$ appropriately depends on the application, but it is generally desirable to minimize these oscillations.
Choosing the $m_k$ that result in an approximant that oscillates 
as little as possible between the sample points yields
\begin{align}
	f(z) = F_0 + \sum_{k=1}^{\floor{M/2}} 
	\left( 
		F_k \e^{\i k z} + F_{M-k} \e^{-\i k z}
	\right)
	+ F_{M/2} \cos \left( M x/2 \right).
	\label{Eq:dft}
\end{align}
The final term vanishes if $M$ is odd.
Moreover, since the trigonometric functions form an orthogonal basis,
truncating \eqref{Eq:dft} at $m<M$ terms provides the $m$-th order trigonometric polynomial approximant that 
is a best fit for the data in the least squares sense.
This corresponds to bandlimiting the signal in the frequency domain with a low-pass filter;
see the discussion in section 10.3.2 of \cite{Sauer2012} for further details.

\begin{figure}[tpb]
	\setlength{\fwidth}{.85\linewidth}
	\setlength{\fheight}{4cm}
	\centering
	\input{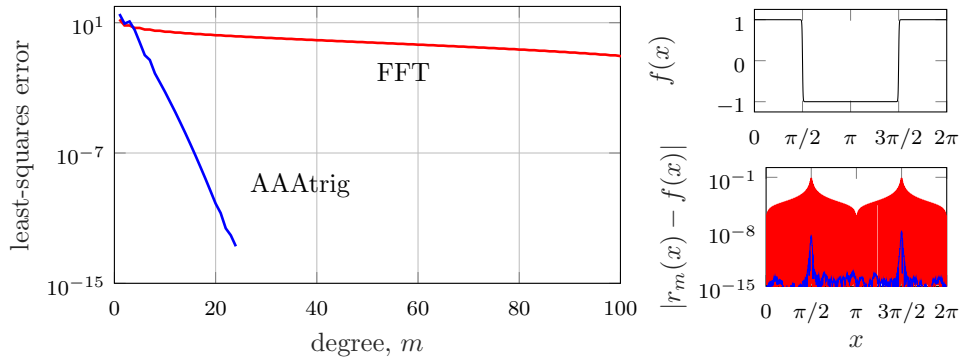}
	\caption{Comparison of AAAtrig (blue) to FFT-based interpolation (red) for $f(x) = \tanh(60 \cos(x)).$
	On the left we compare the least squares error for $m$-th order approximants.
The signal $f$ is plotted at the top right;
on the bottom right we plot the error for the converged AAAtrig approximant and the order 
1,000 trigonometric interpolating polynomial evaluated at a grid 10 times finer than the sample points.}
\label{Fig:dft}
\end{figure}

Figure \ref{Fig:dft} compares the FFT-based approximant \eqref{Eq:dft} to the rational approximant 
computed by AAAtrig.
The function $f(x) = \tanh(60 \cos(x))$ is sampled at 1,000 equi-spaced points in the interval $\left[0,2\pi \right]$;
the function is illustrated in the top right panel of figure \ref{Fig:dft}.
We then approximate $f$ using \eqref{Eq:dft} truncated at different values of $m$ to obtain an
$m$-th order trigonometric polynomial.
The error is compared to that of a type $(m,m)$ trigonometric rational function computed by AAAtrig
in the left panel of figure  \ref{Fig:dft}.
AAAtrig clearly provides an excellent approximant with a type (22,22) trigonometric rational function.
The FFT method does not reach comparable accuracy until nearly 1,000 terms are used.

When an order 1,000 trigonometric polynomial is used, the approximant interpolates the data at each sample point.
The error between this interpolant and $f$ sampled on a grid 10 times finer than the original sampling grid is compared
to that of AAAtrig in the bottom right window of figure \ref{Fig:dft}.
The error of the FFT-based interpolant is large ($10^{-1}$) near the transition points at $\pi/2$ and $3\pi/2$
where there are significant effects of the Gibbs phenomenon.
In contrast, AAAtrig maintains accuracy of $10^{-8}$ even near these points.
If we were to plot the poles and zeros of the AAAtrig approximant then we would see that they cluster exponentially
near $\pi/2$ and $3\pi/2$ where the gradient is almost singular.
This improved accuracy is especially notable when we consider that the AAAtrig approximant only uses 22 terms
whereas the converged FFT interpolant uses 1,000 terms.

This example shows that it is sometimes advantageous to use a trigonometric rational approximant computed 
by AAAtrig as opposed to an trigonometric polynomial approximant computed with an FFT.
AAAtrig is also appropriate when there are discontinuities or singularities in the signal.
Of course, there is no comparison in speed --- the FFT can be computed in $M\log(M)$ flops whereas AAAtrig uses 
approximately $\mathcal{O}(M m^3)$ flops ---
but if the approximant is to be evaluated a large number of times then AAAtrig compression 
may provide a more suitable approximant.
The poles and zeros computed by AAAtrig may also provide useful information for signal classification and diagnostics.
Further applications of AAA to signal processing have been developed in \cite{Wilber2020}.


\subsection{Compression of periodic harmonic functions}
\label{Sec:harmonic}
%
Since the initial discovery of the AAA algorithm, a number of works have applied it to represent harmonic functions, 
especially in domains with corners.
Motivated by Newman's rational approximation of $|x|$ \cite{Newman1964},
Gopal and Trefethen showed that more general corner singularities of the form $x^\alpha$ could be approximated
with root-exponential accuracy by exponentially clustering poles near the corner \cite{Gopal2019a}.
These findings guide a new powerful class of methods to solve Laplace's equation (and other elliptic PDEs) termed
`lightning solvers'.
The strategy is to express the solution as a rational function as the sum of a polynomial (the Runge part) and a series of 
partial fractions (the Newman part).
The poles of the partial fractions are prescribed and cluster exponentially around the corner.
The relevant coefficients are then found by collocating on the boundary and solving a least-squares problem.

These ideas generalise quite naturally to periodic domains --- periodizing the lightning solver of \cite{Gopal2019b} results
in an ansatz of the form
\begin{align}
	f(z) &= \underbrace{\sum^{n_1}_{j=1} a_j \cot\left( \frac{z - z_j}{2}  \right) }_{\textnormal{periodic Newman part}}
	+ \underbrace{\sum^{n_2}_{j=1} b_j \tan\left( \frac{z - z^\ast}{2}  \right)^j}_{\textnormal{periodic Runge part}},
	\label{Eq:lightning}
\end{align}
%
where $n_1$ and $n_2$ represent the number of terms included in each part.
The poles are located at $z=z_j+2 \pi k$ for $k\in \mathbb{Z}$.
The coefficients $a_j$ and $b_j$ are then found by enforcing the boundary condition, which can be of
Dirichlet or Neumann form, and solving a least-squares problem.
When forming the least-squares problem, the Runge part takes Vandermonde structure, which generally leads
to an ill-conditioned matrix.
This issue can be overcome by using the new ``Vandermonde with Arnoldi'' algorithm \cite{Brubeck2019} to
construct a polynomial basis that is discretely orthogonal with respect to the collocation points.
Restructuring the least-squares problem in terms of this orthogonal 
basis drastically improves the condition number of the problem.

\begin{figure}[t]
	\centering
		\includegraphics[width=\linewidth]{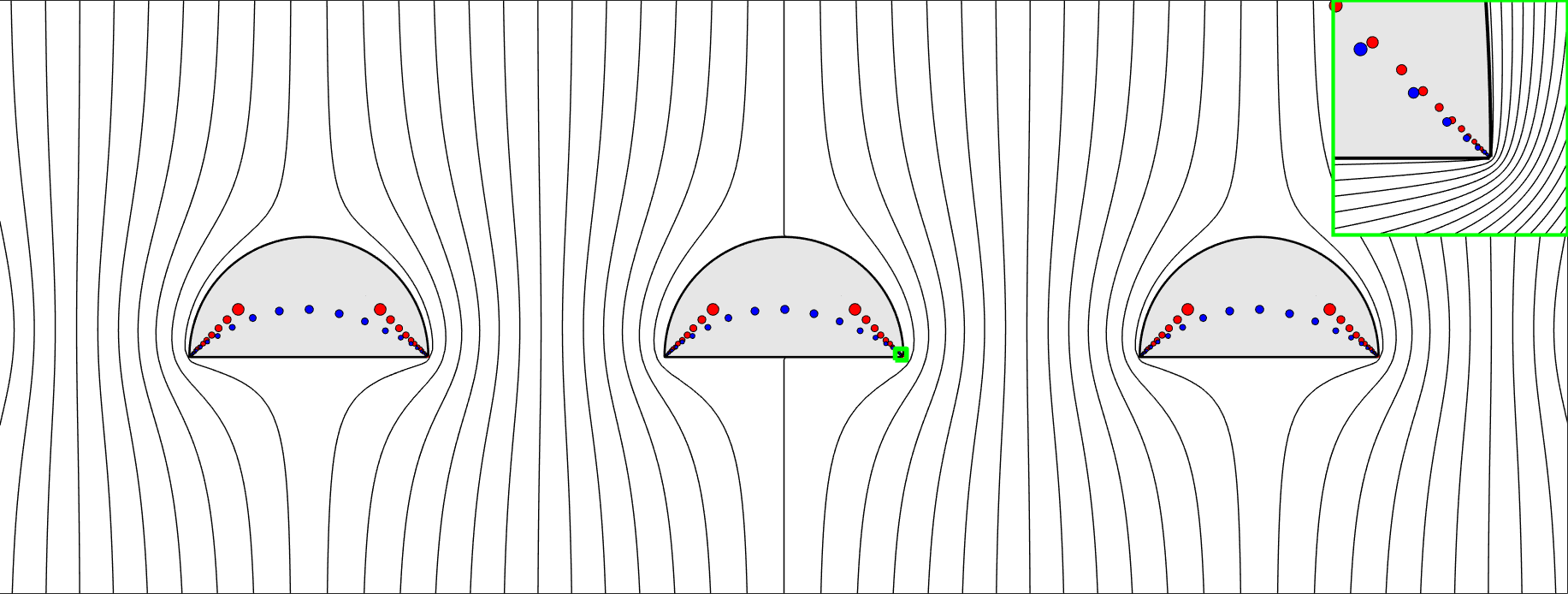}
	\caption{Potential flow through a periodic array of semicircles solved with the periodic lightning
	solver \eqref{Eq:lightning} then compressed with AAAtrig.
	The black lines indicate the streamlines of the flow.
The red circles represent the poles used by the lightning solver and the blue circles represent the poles computed by AAAtrig.
The size of the circles correspond to the size of the residue of each pole.
The close-up view in the inset green box illuminates the exponential clustering of the poles near the corners.
}
\label{Fig:lightning}
\end{figure}
An example of AAAtrig compression is illustrated in figure \ref{Fig:lightning}, which shows the 
potential flow past a periodic array of semicircular obstacles.
Such geometries commonly arise in applications in aerodynamics and microfluidics \cite{Baddoo2020TCFD}.
In this case, the complex function $f(z)$ represents 
the perturbation to a uniform flow induced by the boundaries and
the boundary condition is enforced by requiring that the stream function vanishes on the boundaries, i.e.
$\Im[f(z) + \i z]=0$.
Expressing $f$ in the form \eqref{Eq:lightning} 
and finding the coefficients results in a solution of accuracy $10^{-5}$ using 122 poles.
The solution is then compressed with AAAtrig by evaluating $f(z)$ at 1,000 sample points on the boundary.
The resulting trigonometric rational function has only 36 poles
and is therefore much faster to evaluate than the original lightning solution.
In both cases, the poles cluster exponentially close to the corners to capture 
the (square-root) singularities of the solution there, as illustrated in the green inset box in figure \ref{Fig:lightning}.
Moreover, AAAtrig selects poles that seem to represent a branch cut that connects the two corners.
As we shall see in the next example, AAA and AAAtrig select poles that approximate branch cuts in some sort of optimal way
such that the normal derivative of a potential gradient are balanced \cite{Stahl2012,Trefethen2020Singularity}.
%
%
%
%
%
%

\subsection{Periodic and multi-valued conformal maps}
Conformal maps are another application amenable for AAA compression.
Constructing conformal maps is challenging and even when they are constructed, their evaluation can require 
expensive integrals.
For example, Schwarz--Christoffel (S--C) formulae are a class of formula that represent conformal maps from canonical domains 
(typically the unit disk) to polygons but must be computed using quadrature \cite{Driscoll2002}.
An alternative approach to representing conformal mappings was suggested by Gopal and Trefethen \cite{Gopal2019b}
--- find the boundary correspondence function between the two domains and then approximate the mapping using AAA.
This compression typically results in a speed up of the order of 10--1,000 and has the additional
advantage that both the forward and backward mappings can be computed: the S--C formulae typically only provide the
forward map.

Again, these ideas generalise naturally to periodic domains with minor modifications.
The periodic S--C formula \cite{Baddoo2019c} can be used to determine
the boundary correspondence function between the unit disc (the $\zeta$-plane) and a periodic array of polygons 
(the $z$-plane).
When there is one $N$-sided polygon per period window, the formula takes the form
\begin{align}
	f(\zeta) = A\int^{\zeta} \frac{\prod_{j=1}^{N} \left( \zeta^\prime - \zeta_j \right)^{\beta_j}
	\d \zeta^\prime }
	{\left( \zeta^\prime - a_{\infty^-} \right)
\left( \zeta^\prime - a_{\infty^+} \right)
\left( \zeta^\prime - 1/\overline{a_{\infty^+}} \right)
\left( \zeta^\prime - 1/\overline{a_{\infty^-}} \right)
} + B.
\label{Eq:sc}
\end{align}
In the above, $\beta_j$ are the exterior angles, 
$\zeta_j$ are the pre-images of the corners, $a_{\infty^\pm}$ are the pre-images of $\pm \i\infty$,
$A$ is a scaling and rotation constant and $B$ is a translation constant.
Generally these parameters must be found numerically; 
we calculate them using an adapted version of the SC Toolbox \cite{Driscoll1996}.
As an aside, the periodic S--C formula \eqref{Eq:sc} generalises to an arbitrary number of boundaries per period
window by employing the transcendental Schottky--Klein prime function \cite{Crowdy2007,CrowdyBook}.

\begin{figure}[t]
	\centering
	\begin{tikzpicture}[scale=.95,every node/.style={transform shape}]
	\begin{scope}
	\node at (0,4) {\includegraphics[width=.9\linewidth,trim = .4cm 0 .4cm 0,clip]{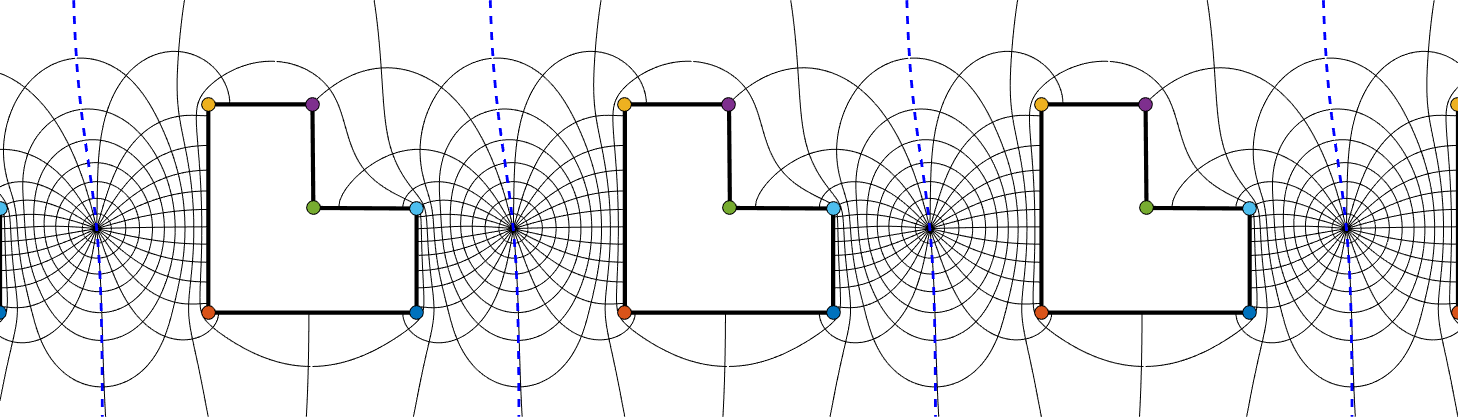}};
	\begin{scope}[shift = {(0,.5)}]
		\node at (-2.25,0) {\includegraphics[width=.3\linewidth]{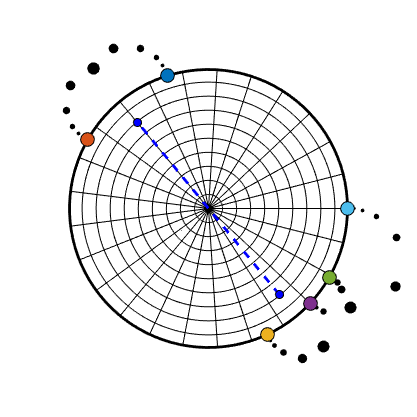}};
		\myPathTextLeftBelow{$z =  f(\zeta)$}{(-5.5,1.5)}{(-4,0)}{-40};
	\end{scope}
	\draw[red!20!white,rounded corners=3mm, line width = 1mm] (-.5\linewidth,-1.5)--(-.2,-1.5)--(-.2,2)--(.5\linewidth,2)%
	--(.5\linewidth,6.5)--(-.5\linewidth,6.5)--cycle;
	\node at (0,6.15) {(a) The forward map};
	\end{scope}
	%
	%
\begin{scope}
	\node at (0,-4) {\includegraphics[width=.9\linewidth,trim = .4cm 0 .4cm 0,clip]{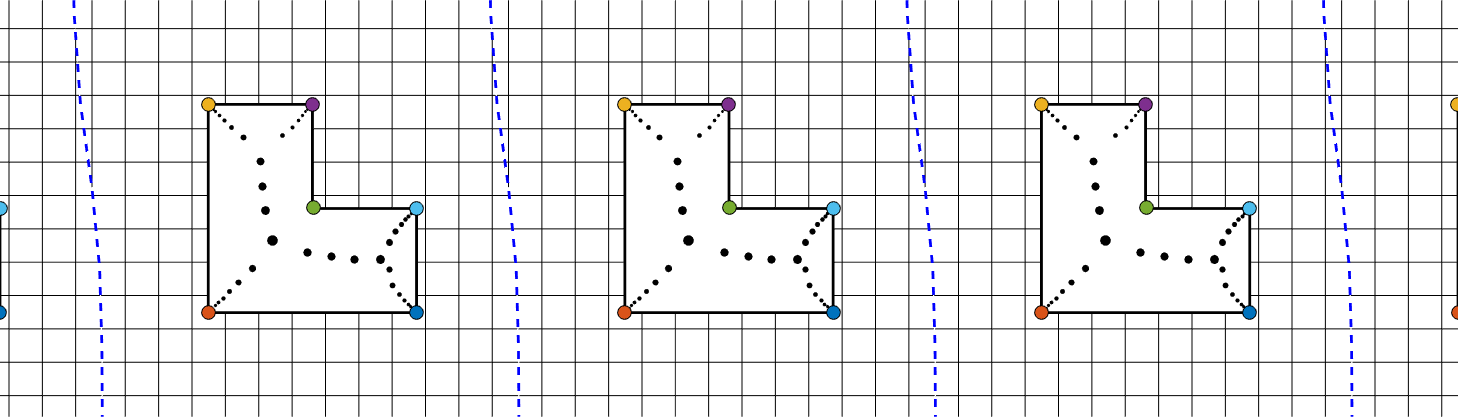}};
\begin{scope}[shift={(0,-.5)}]
		\node at (2.25,0) {\includegraphics[width=.23\linewidth]{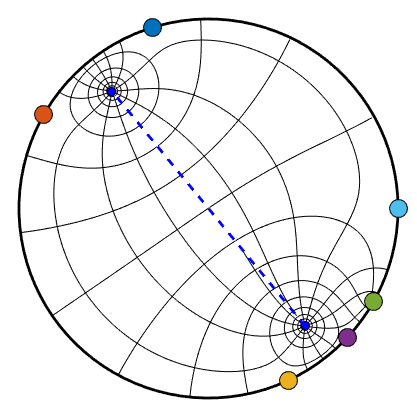}};
		\myPathTextLeftAbove{$\zeta =  f^{-1}(z)$}{(4,0)}{(5.5,-1.5)}{40};
	\end{scope}
	\draw[blue!20!white,rounded corners=3mm, line width = 1mm] (-.5\linewidth,-2)--(.3,-2)--(.3,1.5)--(.5\linewidth,1.5)%
		--(.5\linewidth,-6.5)--(-.5\linewidth,-6.5)--cycle;
	\node at (0,-6.15) {(b) The backward map};
\end{scope}
\end{tikzpicture}
\caption{
	Conformal maps between the unit disc and a periodic array of L-shapnd polygons.
	Figure (a) shows the (multi-valued) forward map from the disc to the periodic domain;
	figure (b) shows the (periodic) backward map from the periodic domain to the disc.
	The forward map is constructed using the periodic Schwarz--Christoffel formula \eqref{Eq:sc} and then compressed
	using AAA \eqref{Eq:branch}.
	The backward map is computed by applying AAAtrig directly to the boundary correspondence function.
	The black dots correspond to the poles computed by AAA and AAAtrig and their size corresponds to the 
	size of the reside at each pole.
	The coloured dots correspond to the corners and pre-images of the corners.
	The dashed blue line represents the branch cut and the endpoints of the branch cut are $a_{\infty^\pm}$.
	A close-up of the structure of the poles in the L-shaped domain is in figure \ref{Fig:tap1}.
}%
\label{Fig:sc}
\end{figure}
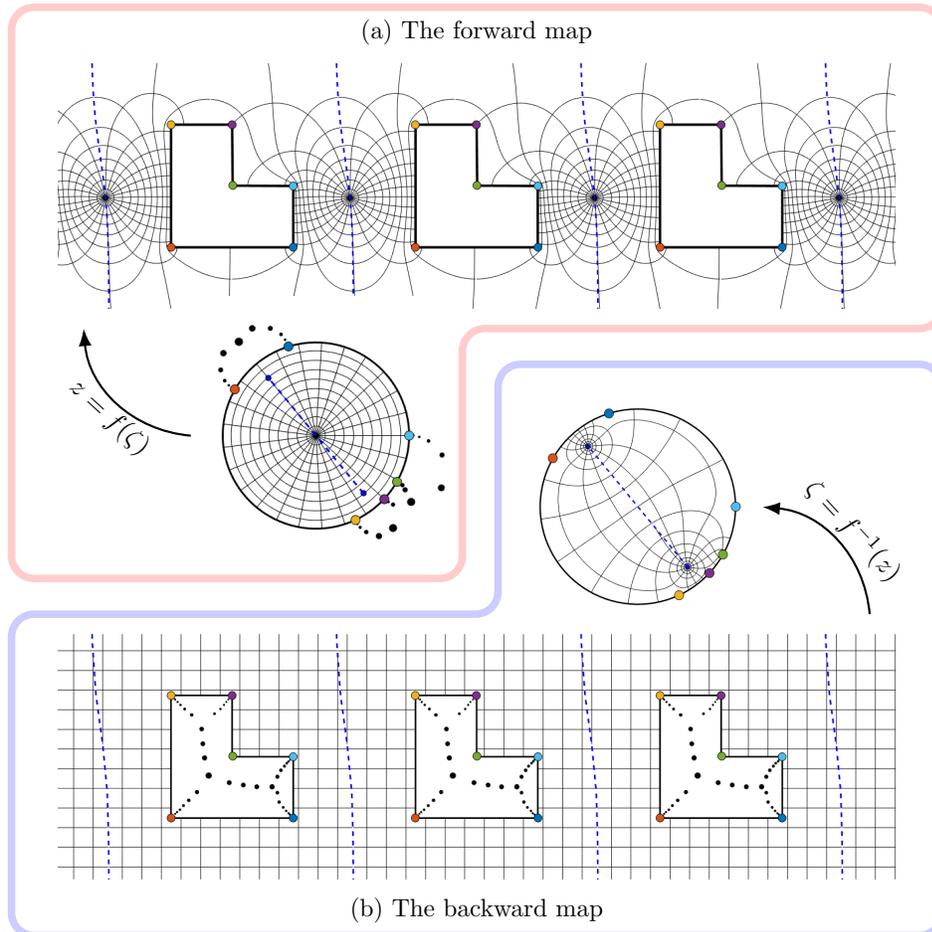
Since the integrand in \eqref{Eq:sc} contains simple poles at $a_{\infty^\pm}$, the resulting $f$ possess a branch cut.
Crossing this branch cut corresponds to traversing between adjacent period windows.
As such, $f$ is multi-valued and cannot be represented with polynomial AAA.
However, since the period of the polygonal domain is known (it is specified in the parameter problem), we know the 
strength of the branch points and can thus extract the branch cut from the conformal map.
As such, when the period is $2 \pi$ we may express the conformal map as
\begin{align}
	z = f(\zeta) = \tilde{f}(\zeta) + \i \log \left( \frac{\zeta - a_{\infty^-}}{\zeta - a_{\infty^+}}  \right)
	\label{Eq:branch}
\end{align}
where $\tilde{f}$ is analytic in the unit disk and can therefore be represented with polynomial AAA.
The result is illustrated in figure \ref{Fig:sc}a --- AAA approximates the branch cuts at the pre-images of the corners
with clusters of poles.
Note that this compression is purely adaptive: the poles and zeros are not user specified and arise naturally
from the AAA algorithm's choice of support points and weights.
Typical tests show that the AAA representation is around 10--100 times 
faster than using quadrature methods to evaluate the integrals.
If AAA is applied directly without first extracting the branch cut then the resulting poles and zeros 
interlace between $a_{\infty^\pm}$ to approximate the branch cut, and the approximation breaks down in a neighbourhood of
the branch points.
Extracting the branch points as in \eqref{Eq:branch} results in a representation that is faithful throughout the entire domain.
%
%

There are typically no constructive formulae for inverse S--C mappings (polygon to disc) 
and they are generally computed with with Newton iteration.
This approach is usually prohibitively expensive in applications.
In our case, the inverse mapping $f^{-1}$ is periodic -- every polygon is effectively mapped to the same disk
(although technically they are on different sheets of an infinite, periodic Riemann surface).
As such, the conformal map can be represented with AAAtrig.
Since $\pm \i \infty$ are mapped to different points $a_{\infty^\pm}$,
it is necessary to use the odd version of AAAtrig \eqref{Eq:cstDef}.
The result is illustrated in figure \ref{Fig:sc}b.
Again, the poles computed by AAAtrig cluster exponentially close to the corners of the polygon.
Moreover, the poles align along some approximate branch cut that connects each of the corners.
We emphasize that the locations of these poles are not specified by the user but are selected
adaptively and automatically by AAAtrig;
we will discuss the locations of these poles in the next section.
Compressing with AAAtrig resulted in speed-up in the range 100--10,000.

Periodic conformal maps between more general non-polygonal shapes can also be constructed by using the lighting solver
of \eqref{Eq:lightning} with appropriately defined boundary values using the approach of \cite{Trefethen2020Conf}.
Again, the number of poles used to represent these maps can be drastically reduced using AAAtrig compression.

%

\begin{figure}[tpb]
	\begin{subfigure}{.35\linewidth}
	\setlength{\fwidth}{.9\linewidth}
	\centering
%
%
\definecolor{mycolor1}{rgb}{0.00000,0.44700,0.74100}%
\definecolor{mycolor2}{rgb}{0.85000,0.32500,0.09800}%
\definecolor{mycolor3}{rgb}{0.92900,0.69400,0.12500}%
\definecolor{mycolor4}{rgb}{0.49400,0.18400,0.55600}%
\definecolor{mycolor5}{rgb}{0.46600,0.67400,0.18800}%
\definecolor{mycolor6}{rgb}{0.30100,0.74500,0.93300}%
\begin{tikzpicture}

\begin{axis}[%
width=\fwidth,
height=\fwidth,
at={(0\fwidth,0\fwidth)},
scale only axis,
xmin=-0.1,
xmax=3.24159265358979,
ymin=-0.1,
ymax=3.24159265358979,
axis line style={draw=none},
ticks=none,
axis x line*=bottom,
axis y line*=left
]

\addplot[area legend, line width=1pt, draw=black, fill=white, forget plot]
table[row sep=crcr] {%
x	y\\
3.14159265358979	0\\
0	0\\
0	3.14159265358979\\
1.5707963267949	3.14159265358979\\
1.5707963267949	1.5707963267949\\
3.14159265358979	1.5707963267949\\
3.14159265358979	0\\
}--cycle;
\addplot[scatter, only marks, mark=*, color=mycolor1, mark options={}, scatter/use mapped color={mark options={}, draw=black, fill=mycolor1}, visualization depends on={\thisrow{size} \as \perpointmarksize}, scatter/@pre marker code/.append style={/tikz/mark size=\perpointmarksize}] table[row sep=crcr]{%
x	y	size\\
3.14151918840798	5.61167107797644e-05	0.00333326810912449\\
3.14133458641814	0.000242739191954355	0.00913005991421688\\
3.1410101564642	0.000561031293715243	0.0169555295161315\\
3.14051660035596	0.00105525955062988	0.0270439205742125\\
3.13978581161179	0.00180112862010407	0.0401343661836084\\
3.13868405662396	0.00293090463088746	0.0575907813741018\\
3.13698024167596	0.0046555430972076	0.0814348354209758\\
3.1343435299771	0.00728471681785152	0.114118229639666\\
3.13033497242917	0.011283506380728	0.158599409548766\\
3.12433686566335	0.0172995776606419	0.218502111784393\\
3.11543025148524	0.0262483503595915	0.298630727532072\\
3.10230085740451	0.0394326200033391	0.405183140024539\\
3.08311042158771	0.0587049324434441	0.546046189354506\\
3.05533764683204	0.0867143161343548	0.731217742151455\\
3.01559813821953	0.127305617124384	0.973418697873154\\
2.95958273794778	0.186317491528042	1.28904048726504\\
2.88243587175003	0.273529929167662	1.70036118324152\\
2.78112389434835	0.409654353693844	2.24466655462553\\
2.72986825341797	0.653002483758075	2.91174281992298\\
};
\addplot[scatter, only marks, mark=*, color=mycolor2, mark options={}, scatter/use mapped color={mark options={}, draw=black, fill=mycolor2}, visualization depends on={\thisrow{size} \as \perpointmarksize}, scatter/@pre marker code/.append style={/tikz/mark size=\perpointmarksize}] table[row sep=crcr]{%
x	y	size\\
6.72254547279518e-05	8.03348346242507e-05	0.00366080202724357\\
0.000291310862678085	0.00029239362018676	0.0102379755973844\\
0.000669357132039421	0.000663789893019033	0.0190208593483989\\
0.00124891664176453	0.00123625204284444	0.0303451921539036\\
0.00211434297792684	0.00209427408237157	0.0450481561622722\\
0.00342184211527624	0.00340017626174967	0.0647150162384277\\
0.00543934530857672	0.00542454937128984	0.0917398910920282\\
0.00855679766167964	0.00854006108546088	0.128900961511651\\
0.0133097691491949	0.0132701808760955	0.179467056332008\\
0.0204804521477783	0.0203900484496526	0.247814676214309\\
0.0311897045717296	0.0310059681364339	0.339539501384194\\
0.0470306905599869	0.0466877049475637	0.461784735435939\\
0.0702580719894882	0.0696515384471362	0.623673770244596\\
0.104031102264803	0.103003040445026	0.836766438853094\\
0.152729830459996	0.151077240847064	1.1156346274351\\
0.222367529250206	0.21993419831074	1.47864072704709\\
0.321148133835012	0.3181924355818	1.94926747562916\\
0.460578921886198	0.458996617219852	2.56011690964305\\
0.659829771200415	0.668056404027794	3.3724509814643\\
};
\addplot[scatter, only marks, mark=*, color=mycolor3, mark options={}, scatter/use mapped color={mark options={}, draw=black, fill=mycolor3}, visualization depends on={\thisrow{size} \as \perpointmarksize}, scatter/@pre marker code/.append style={/tikz/mark size=\perpointmarksize}] table[row sep=crcr]{%
x	y	size\\
0.000104958744425792	3.14146798361863	0.00509957556612071\\
0.000455410947099948	3.14113742495082	0.0142915278097433\\
0.00104521722899053	3.14055790238864	0.0265531127166383\\
0.00194693220248203	3.13966279151748	0.0423571257395972\\
0.00328987789834748	3.13832100040591	0.0628526954683824\\
0.00532175830923139	3.13628375946606	0.0902586713107206\\
0.00847666812577765	3.13312260339444	0.128051004942086\\
0.0133556067707259	3.12824761491049	0.180077769865621\\
0.0208088511159218	3.12083771439758	0.250960511564363\\
0.0320582303377069	3.10967749020602	0.346792035870633\\
0.0488673733846595	3.0930321367061	0.475425911060207\\
0.0737651672338798	3.06844783521748	0.646934257439139\\
0.110355250714095	3.03248108897812	0.874179533092439\\
0.163770448419164	2.98035285768254	1.17354951771794\\
0.241438352897586	2.90545099633268	1.56629439801975\\
0.354903604907529	2.7981921533024	2.08279435449538\\
0.526983371045466	2.6397727301256	2.78557860665771\\
};
\addplot[scatter, only marks, mark=*, color=mycolor4, mark options={}, scatter/use mapped color={mark options={}, draw=black, fill=mycolor4}, visualization depends on={\thisrow{size} \as \perpointmarksize}, scatter/@pre marker code/.append style={/tikz/mark size=\perpointmarksize}] table[row sep=crcr]{%
x	y	size\\
1.57053655442014	3.14137445546565	0.00883785935215678\\
1.56978381512831	3.14060644651538	0.0257727480017938\\
1.56849088306968	3.13931621901457	0.0480129011277357\\
1.56651403089216	3.13732877818099	0.0766298914537845\\
1.5635687369211	3.13435063170802	0.113739569191714\\
1.55911111071918	3.12984723236579	0.163270570165817\\
1.55221201569999	3.12290846949281	0.231247113659416\\
1.54160087265468	3.11223189041884	0.324526238852066\\
1.52555136299736	3.09596114284173	0.451242823183731\\
1.50161076578646	3.07125046920631	0.622408308547034\\
1.46663599805695	3.03380598417285	0.851700455371326\\
1.41779673755547	2.97704311745233	1.15383115148735\\
1.35704625607735	2.89464503739941	1.52747550077572\\
1.27647553110353	2.79950074479159	1.94665031918176\\
1.12155201674761	2.67028363912561	2.56270732026551\\
};
\addplot[only marks, mark=*, mark options={}, mark size=1.5000pt, draw=black, fill=mycolor5] table[row sep=crcr]{%
x	y\\
};
\addplot[scatter, only marks, mark=*, color=mycolor6, mark options={}, scatter/use mapped color={mark options={}, draw=black, fill=mycolor6}, visualization depends on={\thisrow{size} \as \perpointmarksize}, scatter/@pre marker code/.append style={/tikz/mark size=\perpointmarksize}] table[row sep=crcr]{%
x	y	size\\
3.1415680756113	1.57077241333153	0.00158433918140637\\
3.14148594577958	1.57069130142324	0.00478542742020952\\
3.14134770318853	1.57055431313936	0.00893711530748732\\
3.14113652386677	1.57034319684199	0.0142753406068936\\
3.14082199859876	1.57002668700233	0.0211968688308828\\
3.1403459258363	1.56954701524745	0.0304444054431076\\
3.13960977786447	1.56880685081214	0.0431375005553916\\
3.13848025502182	1.56767244622572	0.0605011277099004\\
3.1367633615654	1.56595046688502	0.0841032117387513\\
3.13417785928476	1.56335788650162	0.115993846731611\\
3.13032063167701	1.55949361851827	0.158776412631766\\
3.1246183176677	1.55378767331277	0.215781947923876\\
3.11626058918805	1.54543263029436	0.291271207292344\\
3.1041165642529	1.53328946992088	0.390652547731744\\
3.08662981623034	1.51575871391948	0.520732496040597\\
3.06169532539777	1.49058732469064	0.690042361649232\\
3.02655304288038	1.45456310854381	0.909214923915237\\
2.97777508404561	1.40302223903809	1.19136445754441\\
2.91153828973842	1.3288531399574	1.55284054159472\\
2.8252885686994	1.22009915614282	2.01417030414574\\
2.72596551590272	1.05638879989009	2.59280606937728\\
};
\addplot[only marks, mark=*, mark options={}, mark size=3.0619pt, draw=black, fill=white] table[row sep=crcr]{%
x	y\\
2.21256315540455	0.797408763298954\\
1.86115191730208	0.840278907610228\\
1.49635919552045	0.906223435783683\\
2.604233532583	0.803612715059937\\
0.966000614964507	1.08231798389465\\
0.870440911301789	1.53500317465281\\
0.822472038344663	1.90002008079039\\
0.788519192796664	2.2849450163028\\
};
\end{axis}

\begin{axis}[%
width=\fwidth,
height=\fwidth,
at={(0\fwidth,0\fwidth)},
scale only axis,
xmin=0,
xmax=1,
ymin=0,
ymax=1,
axis line style={draw=none},
ticks=none,
axis x line*=bottom,
axis y line*=left
]
\end{axis}
\end{tikzpicture}%
	\vspace{.7cm}
	\caption{}
	\label{Fig:tap1}
\end{subfigure}
\hfill
\begin{subfigure}{.6\linewidth}
	\setlength{\fwidth}{.7\linewidth}
	\setlength{\fheight}{6cm}
	\centering
%
%
\definecolor{mycolor1}{rgb}{0.00000,0.44700,0.74100}%
\definecolor{mycolor2}{rgb}{0.85000,0.32500,0.09800}%
\definecolor{mycolor3}{rgb}{0.92900,0.69400,0.12500}%
\definecolor{mycolor4}{rgb}{0.49400,0.18400,0.55600}%
\definecolor{mycolor5}{rgb}{0.46600,0.67400,0.18800}%
\definecolor{mycolor6}{rgb}{0.30100,0.74500,0.93300}%
\begin{tikzpicture}[%
trim axis left, trim axis right
]

\begin{axis}[%
width=0.982\fwidth,
height=\fheight,
at={(0\fwidth,0\fheight)},
scale only axis,
xmin=1,
xmax=4.5,
xlabel style={font=\color{white!15!black}},
xlabel={$\sqrt{k}$},
ymode=log,
ymin=0.0001,
ymax=1,
yminorticks=true,
ylabel style={font=\color{white!15!black}},
ylabel={$d_k$},
axis background/.style={fill=white},
xmajorgrids,
ymajorgrids
]
\addplot [color=mycolor1, line width=1.0pt, forget plot]
  table[row sep=crcr]{%
1	0.000354289402036691\\
1.41421356237309	0.000808739150683667\\
2	0.00255122367494816\\
2.23606797749979	0.00412918130686782\\
2.44948974278318	0.00655350480222622\\
2.64575131106459	0.0102770079433858\\
2.82842712474619	0.0159390370273159\\
3	0.024434352952331\\
3.16227766016838	0.0370600483065019\\
3.3166247903554	0.0556666575949549\\
3.46410161513775	0.0828639882767816\\
3.60555127546399	0.122308212371163\\
3.74165738677394	0.179112640689614\\
3.87298334620742	0.260464617638012\\
4	0.376803476263166\\
4.12310562561766	0.545668778554707\\
4.24264068711928	0.771964523466633\\
};
\addplot[only marks, mark=*, mark options={}, mark size=1.5811pt, draw=mycolor1, fill=mycolor1] table[row sep=crcr]{%
x	y\\
1	0.000354289402036692\\
1.4142135623731	0.000808739150683667\\
1.73205080756888	0.00150713744603673\\
2	0.00255122367494816\\
2.23606797749979	0.00412918130686782\\
2.44948974278318	0.00655350480222623\\
2.64575131106459	0.0102770079433858\\
2.82842712474619	0.0159390370273159\\
3	0.024434352952331\\
3.16227766016838	0.0370600483065019\\
3.3166247903554	0.0556666575949549\\
3.46410161513775	0.0828639882767816\\
3.60555127546399	0.122308212371163\\
3.74165738677394	0.179112640689614\\
3.87298334620742	0.260464617638012\\
4	0.376803476263167\\
4.12310562561766	0.545668778554707\\
4.24264068711928	0.771964523466633\\
};
\addplot [color=mycolor2, line width=1.0pt, forget plot]
  table[row sep=crcr]{%
1	0.000412742108150076\\
1.41421356237309	0.000942685521415417\\
2	0.00297597549727847\\
2.23606797749979	0.00482391978300262\\
2.44948974278318	0.00768194072272733\\
2.64575131106459	0.0120893105496686\\
2.82842712474619	0.0187948837530098\\
3	0.0288998788225248\\
3.16227766016838	0.0439791738365936\\
3.3166247903554	0.0662693567783783\\
3.46410161513775	0.0989319639334678\\
3.60555127546399	0.146397051129282\\
3.74165738677394	0.214827218513625\\
3.87298334620742	0.31275928387728\\
4	0.452086883162063\\
4.12310562561766	0.650239062118787\\
4.24264068711928	0.938975338294324\\
};
\addplot[only marks, mark=*, mark options={}, mark size=1.5811pt, draw=mycolor2, fill=mycolor2] table[row sep=crcr]{%
x	y\\
1	0.000412742108150076\\
1.4142135623731	0.000942685521415418\\
1.73205080756888	0.00175730244736455\\
2	0.00297597549727847\\
2.23606797749979	0.00482391978300262\\
2.44948974278318	0.00768194072272733\\
2.64575131106459	0.0120893105496686\\
2.82842712474619	0.0187948837530098\\
3	0.0288998788225248\\
3.16227766016838	0.0439791738365937\\
3.3166247903554	0.0662693567783783\\
3.46410161513775	0.0989319639334679\\
3.60555127546399	0.146397051129282\\
3.74165738677394	0.214827218513625\\
3.87298334620742	0.31275928387728\\
4	0.452086883162063\\
4.12310562561766	0.650239062118787\\
4.24264068711928	0.938975338294323\\
};
\addplot [color=mycolor3, line width=1.0pt, forget plot]
  table[row sep=crcr]{%
1	0.000643919439431223\\
1.41421356237309	0.00147077840073609\\
2	0.00463972101980897\\
2.23606797749979	0.00751701192754967\\
2.44948974278318	0.0119831403574484\\
2.64575131106459	0.018880208938633\\
2.82842712474619	0.0293900627021603\\
3	0.0452361336344989\\
3.16227766016838	0.0688922635825469\\
3.3166247903554	0.103881973179953\\
3.46410161513775	0.155188965110852\\
3.60555127546399	0.229823914246594\\
3.74165738677394	0.337720832259478\\
3.87298334620742	0.493842558285528\\
4	0.727691355551273\\
};
\addplot[only marks, mark=*, mark options={}, mark size=1.5811pt, draw=mycolor3, fill=mycolor3] table[row sep=crcr]{%
x	y\\
1	0.000643919439431223\\
1.4142135623731	0.00147077840073609\\
1.73205080756888	0.00274133409478453\\
2	0.00463972101980897\\
2.23606797749979	0.00751701192754967\\
2.44948974278318	0.0119831403574484\\
2.64575131106459	0.018880208938633\\
2.82842712474619	0.0293900627021603\\
3	0.045236133634499\\
3.16227766016838	0.0688922635825469\\
3.3166247903554	0.103881973179953\\
3.46410161513775	0.155188965110852\\
3.60555127546399	0.229823914246594\\
3.74165738677394	0.337720832259478\\
3.87298334620742	0.493842558285527\\
4	0.727691355551274\\
};
\addplot [color=mycolor4, line width=1.0pt, forget plot]
  table[row sep=crcr]{%
1	0.00141342996593523\\
1.4142135623731	0.00323994215155377\\
2	0.0102315656827259\\
2.23606797749979	0.0165680172158747\\
2.44948974278318	0.0263529003003887\\
2.64575131106459	0.0414056633376877\\
2.82842712474619	0.0642599527092146\\
3	0.0986644047055104\\
3.16227766016838	0.149891094421234\\
3.3166247903554	0.22468961736246\\
3.46410161513775	0.326607130776212\\
3.60555127546399	0.451277746893859\\
3.74165738677394	0.651116454426598\\
};
\addplot[only marks, mark=*, mark options={}, mark size=1.5811pt, draw=mycolor4, fill=mycolor4] table[row sep=crcr]{%
x	y\\
1	0.00141342996593523\\
1.4142135623731	0.00323994215155377\\
1.73205080756888	0.00604306972493447\\
2	0.0102315656827259\\
2.23606797749979	0.0165680172158747\\
2.44948974278318	0.0263529003003887\\
2.64575131106459	0.0414056633376878\\
2.82842712474619	0.0642599527092147\\
3	0.0986644047055104\\
3.16227766016838	0.149891094421234\\
3.3166247903554	0.22468961736246\\
3.46410161513775	0.326607130776212\\
3.60555127546399	0.451277746893859\\
3.74165738677394	0.651116454426597\\
};
\addplot[only marks, mark=*, mark options={}, mark size=1.5811pt, draw=mycolor5, fill=mycolor5] table[row sep=crcr]{%
x	y\\
};
\addplot [color=mycolor6, line width=1.0pt, forget plot]
  table[row sep=crcr]{%
1	0.00014972269518251\\
1.41421356237309	0.000344341848376428\\
2	0.00108915312307594\\
2.23606797749979	0.00176496726142833\\
2.44948974278318	0.00280888070732859\\
2.64575131106459	0.00440972274144623\\
2.82842712474619	0.00684137557242812\\
3	0.0105028362161244\\
3.16227766016838	0.0159627595480505\\
3.3166247903554	0.0240296145052495\\
3.46410161513775	0.0358473232895512\\
3.60555127546399	0.053020954202601\\
3.74165738677394	0.0777820822659495\\
3.87298334620742	0.11321248637312\\
4	0.163536763625724\\
4.12310562561766	0.234487399691307\\
4.24264068711928	0.333858526900455\\
4.47213595499958	0.661332761413094\\
};
\addplot[only marks, mark=*, mark options={}, mark size=1.5811pt, draw=mycolor6, fill=mycolor6] table[row sep=crcr]{%
x	y\\
1	0.000149722695182509\\
1.4142135623731	0.000344341848376428\\
1.73205080756888	0.00064294718168617\\
2	0.00108915312307594\\
2.23606797749979	0.00176496726142833\\
2.44948974278318	0.00280888070732859\\
2.64575131106459	0.00440972274144623\\
2.82842712474619	0.00684137557242813\\
3	0.0105028362161244\\
3.16227766016838	0.0159627595480505\\
3.3166247903554	0.0240296145052495\\
3.46410161513775	0.0358473232895512\\
3.60555127546399	0.053020954202601\\
3.74165738677394	0.0777820822659495\\
3.87298334620742	0.11321248637312\\
4	0.163536763625724\\
4.12310562561766	0.234487399691307\\
4.24264068711928	0.333858526900456\\
4.35889894354067	0.47226769910901\\
4.47213595499958	0.661332761413095\\
};
\end{axis}
\end{tikzpicture}%
	\caption{}
	\label{Fig:tap2}
\end{subfigure}
	\caption{
		(a) A close-up of the poles discovered by AAAtrig for the backward mapping in figure \ref{Fig:sc}.
		(b) The distance $d_k$ from each corner to its $k$-th nearest poles on a square-root scale.
		The line color corresponds to the locations of the poles in figure (a).
		Note that the log-distance scales approximately linearly on the $\sqrt{k}$ scale for small $k$.
		Poles with white fill in figure (a) are not directly associated to a particular corner and
		are therefore not plotted in figure (b).
		The reentrant corner does not have an associated cluster of poles.
	}
	\label{Fig:tap}
\end{figure}
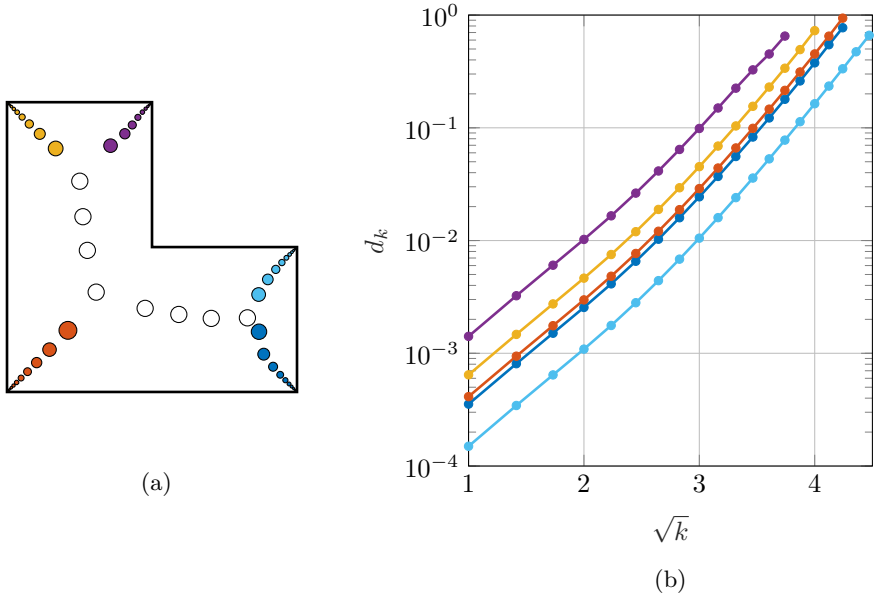

\subsection{Clustering of poles and zeros}
We end with a brief discussion of the distribution of poles selected by AAAtrig.
A recent study \cite{Trefethen2020Singularity} has shown that the optimal distribution of poles and zeros
near a singularity is given by the tapered distribution:
\begin{align}
	d_k \approx \beta \exp( \sigma \sqrt{k}),
	\label{Eq:tape}
\end{align}
where $d_k$ is the distance of the $k$-th pole from the singularity and $\beta>0$ and $\sigma<0$ are constants.
With this distribution, the density of poles on a log scale decays linearly as the distance to the singularity decreases.
In the conformal mapping example we see that AAAtrig selects poles and zeros that obey this tapering law.
Figure \ref{Fig:tap} shows the distributions of poles nearest to five corners of the L-shaped polygon
in the backward map in figure \ref{Fig:sc}.
Again, this distribution of poles was not specified but arose through the application of AAAtrig.
The reentrant corner (green) does not have an associated distribution of poles and is thus excluded.
The log-distance is seen to vary linearly with respect to $\sqrt{k}$ in figure \ref{Fig:tap2}, 
showing that these poles obey the tapering law \eqref{Eq:tape}.

 \section{Conclusions} \label{Sec:conclusions}
We have presented an extension of the AAA algorithm \cite{Nakatsukasa2018} to compute 
rational approximants for periodic functions, implemented in Chebfun.
The sample points need not be equispaced and are not constrained to a particular subset of the complex plane.
As with the original AAA algorithm, AAAtrig exploits a (trigonometric) barycentric representation of
the approximant and uses adaptive, greedy selection of support points.
Thus, AAAtrig achieves the same high levels of robustness, flexibility and accuracy seen in the original AAA algorithm.
As one would expect, AAAtrig is an improvement over AAA when the underlying function is periodic.
We have demonstrated that the periodic approximant is essential is applications to PDEs and conformal maps 
when the approximant must be periodic throughout the entire complex plane.%
{We have also developed a minimax (AAA--Lawson) \cite{Nakatsukasa2019} version of AAAtrig
though we did not report any of the details in this paper in the interests of space.}

Much research is currently dedicated to understanding the theoretical underpinnings of the AAA framework
and we expect that these advances will translate naturally to AAAtrig.
Other lines of enquiry should develop applications to PDEs in the spirit of the new `lightning solvers'.
We have seen how AAAtrig can be used to compress solutions of Laplace's equation, and there is significant scope
to develop solvers for other PDEs in periodic domains, such as the Helmholtz, convected Helmholtz and biharmonic equations.
Signal processing is another promising application, 
and we showed that AAAtrig has some advantages over FFT-based interpolation in section \ref{Sec:fft}.
A parallel study has used a similar version of AAAtrig to perform signal reconstruction for noisy or incomplete data
by combining a similar version of AAAtrig with Prony's method \cite{Wilber2020}.





 \section*{Acknowledgements}%
\label{sec:acknowledgements}
The author acknowledges insightful conversations with Profs L. N. Trefethen and Y. Nakatsukasa
at the ``Complex analysis: techniques, applications and computations'' programme at the 
Isaac Newton Institute for Mathematical Sciences, Cambridge, supported by EPSRC grant no \mbox{EP/R014604/1}.
The author also acknowledges an EPSRC Doctoral Prize Fellowship from Imperial College London.

 \bibliographystyle{siamplain}
 \bibliography{../../../bibliography/library}
\end{document}